\documentclass[11pt, oneside, a4paper]{article}
\usepackage{geometry}
\title{A parallelogram height inequality for Drinfeld modules}
\author{{Liam {\sc Baker}} \and {Richard {\sc Griffon}} \and {Fabien {\sc Pazuki}}}
\date{} 
\usepackage[english]{babel}
\usepackage{csquotes}
\usepackage{lmodern}
\usepackage{amssymb,amscd,amsthm}
\usepackage{mathrsfs}
\usepackage{tikz-cd}
\usepackage{mathtools}
\mathtoolsset{centercolon}
\usepackage{enumerate}
\usepackage{graphicx}
\usepackage{titlesec}
\titleformat{\subsection}[runin]{\bfseries}{\thesubsection.}{0.2em}{}[.\hspace{0.4em}-- ]
\titleformat{\section}{\center\Large\bfseries}{\thesection.}{0.2em}{}[]
\usepackage{url}
\definecolor{BrickRed}{RGB}{153, 0, 0}
\definecolor{TealBlue}{RGB}{84, 181, 183}
\definecolor{Orchid}{RGB}{164,117,181}
\definecolor{Violet}{RGB}{123,59,143}
\usepackage[colorlinks=true, citecolor= Orchid, linkcolor= Violet, urlcolor=TealBlue]{hyperref}
\usepackage{scalerel}
\usepackage[
  backend=biber,
  maxbibnames=99,
  maxcitenames=2
]{biblatex}
\newtheorem{itheo}{Theorem}

\newtheorem*{itheono}{Theorem}
 
\newtheorem{prop}{Proposition}[section]

\newtheorem{lem}[prop]{Lemma}
\newtheorem{theorem}[prop]{Theorem}

{   \theoremstyle{definition}
    \newtheorem{definition}[prop]{Definition}
    \newtheorem{rem}[prop]{Remark}
    \newtheorem{ex}[prop]{Example}}
\newcommand{\F}{\ensuremath{\mathbb{F}}}
\newcommand{\Z}{\ensuremath{\mathbb{Z}}}
\newcommand{\Q}{\ensuremath{\mathbb{Q}}}
\newcommand{\R}{\ensuremath{\mathbb{R}}}
\newcommand{\C}{\ensuremath{\mathbb{C}}}

\newcommand{\KERS}{\mathtt{SG}}
\newcommand{\hTag}{\ensuremath{h_{\mathrm{Tag}, K}}}
\newcommand{\hG}{\ensuremath{h_{\mathrm{Gr}, K}}}

\newcommand{\hdiff}{\ensuremath{\mathrm{h}_{\mathrm{diff}}}}

\newcommand{\covol}{\mathcal{D}}

\DeclareMathOperator{\ord}{ord}
\DeclareMathOperator{\Gal}{Gal}
\newcommand{\dd}{\,\mathrm{d}}
\newcommand{\ee}{\mathrm{e}}
\newcommand{\into}{\hookrightarrow}
\newcommand{\ie}{{i.e.}}

\DeclarePairedDelimiter{\card}{|}{|}
\newcommand{\V}{\mathbf{V}}
\newcommand{\N}{\mathbf{N}}
\renewcommand{\L}{\mathfrak{L}}

\def\p[#1]_#2{
	\setbox0=\hbox{$\scriptstyle{#2}$}
	\setbox2=\hbox{$\displaystyle{#1}$}
	\setbox4=\hbox{${}'\mathsurround=0pt$}
	\dimen0=.5\wd0 \advance\dimen0 by-.5\wd2
	\ifdim\dimen0>0pt
	\ifdim\dimen0>\wd4 \kern\wd4 \else\kern\dimen0\fi\fi
	\mathop{{#1}'}_{\kern-\wd4 #2}}
\def\prodp_#1{\p[\prod]_{#1}}

\begin{document}
\pagestyle{plain}
    \maketitle

    \centerline{\rule{7cm}{0.5pt}}
	\paragraph{Abstract --}%
    We prove an inequality relating the Taguchi heights of four Drinfeld modules arranged in a ``parallelogram of isogenies''.
    This inequality is the analogue for Drinfeld modules of the parallelogram inequality of R\'emond~\cite{Re22} for abelian varieties over number fields and of Griffon--Le~Fourn--Pazuki~\cite{GLP26} for abelian varieties over function fields.
    We also prove a similar inequality for the local graded height at a finite place.
	
	\medskip
	
    
	\smallskip
	\noindent\textit{2020 Math. Subj. Classification:}  
    {Primary 11G09, 11G50, 14G40, 14K02; Secondary 11J93.}

    \centerline{\rule{7cm}{0.5pt}}
\section*{Introduction}

Fix a finite field $\F_q$, and let $K_0=\F_q(T)$ denote the rational function field over $\F_q$.
In this article, we consider Drinfeld $\F_q[T]$-modules $\phi$  defined over a finite field extension $K$ of $K_0$.
There are several notions of heights attached to these objects: our main focus will be on the Taguchi height $\hTag(\cdot)$, we will also prove related properties of the graded height $\hG(\cdot)$. We give precise definitions of these two heights in section \ref{sec:heights} below.
This choice is motivated by the growing relevance of these two heights in Diophantine problems related to Drinfeld modules. 
See, for instance, \cite{Wei20} where Wei gives a formula ``\`a la Colmez'' for the stable Taguchi height of a  CM Drinfeld module over~$K$.
Moreover, in a recent paper \cite{BPR21} Breuer, Pazuki, and Razafinjatovo gave bounds on the variation of the graded height by isogeny in this context. 
\bigskip

Let $\phi$ be a Drinfeld module of rank $r\geqslant 1$ defined  over $K$. To any isogeny $f:\phi\to\phi'$ of Drinfeld modules, one associates a finite dimensional $\F_q$-vector space $\ker(f)\subset \overline{K}$ called the kernel of $f$. 
We can then think of the target Drinfeld module $\phi'$ as ``the quotient'' $\phi/\ker(f)$. 
We let $\KERS_K(\phi)$ denote the set of all kernels of isogenies $f:\phi\to\phi'$ between $\phi$ and other Drinfeld modules $\phi'$ defined over $K$. We refer to section~\ref{sec:Drinfeld_isogenies} for precise definitions.
With this notation at hand, we can state our main result:

\begin{itheo}\label{itheo:main}
    Let $K$ be a finite extension of $K_0 = \F_q(T)$.
    Let $\phi:A\to K\{\tau\}$ be a Drinfeld module of rank $r\geqslant 1$ defined over $K$, with stable reduction at all finite places of $K$. For any two $G, H \in \KERS_K(\phi)$, we have
    \[ \hTag\big(\phi/(G\cap H)\big) + \hTag\big(\phi/(G+H)\big)
        \leqslant \hTag(\phi/G) + \hTag(\phi/H).\]
\end{itheo}
 
This relation is called a ``parallelogram inequality'' because of the shape of the diagram of isogenous Drinfeld modules that one constructs from the data $\phi, G, H$ of the theorem:
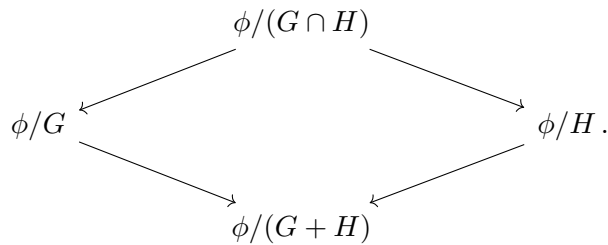
\begin{center}
\begin{tikzcd}
    &  & \phi/(G\cap H) \arrow[lld ] \arrow[rrd ]    &  & \\
    \phi/G \arrow[rrd] &  &  &  & \phi/H\,. \arrow[lld] \\
    &  & \phi/(G+H) &  &                   
\end{tikzcd}
\end{center}
Arrows in this picture represent isogenies of Drinfeld modules.
We underline the fact that Theorem~\ref{itheo:main} involves no extraneous terms (e.g. a dependence on the size of $G, H$). 

It is interesting to note the perfect analogy between Theorem~\ref{itheo:main} and a very recent result of \cite{GLP26} in the context of abelian varieties defined over the function field $K$:
\begin{itheono}[Theorem 7.1 in \cite{GLP26}]
Let $B$ be a semi-stable abelian variety defined over $K$. 
For any finite subgroup schemes $G,H$ of $B$, we have
   \[\hdiff\big(B/(G\cap H)\big) + \hdiff\big(B/(G+H)\big)
        \leqslant \hdiff(B/G) + \hdiff(B/H),\]
where $\hdiff(\cdot)$ denotes the differential height over $K$.
\end{itheono}
This result was itself inspired by the parallelogram inequality of R\'emond (Th\'eor\`eme 1.1 in~\cite{Re22}) in the setting of abelian varieties over a number field, where the role of the differential height is played by the Faltings height. 
The proofs of the respective parallelogram inequalities in these three contexts are based on different strategies.
With Theorem \ref{itheo:main}, the classical Diophantine triptych ``abelian varieties over number fields'' -- ``abelian varieties over function fields'' -- ``Drinfeld modules'' is thus complete as regards the parallelogram inequality. 
\bigskip

Both the graded and the Taguchi heights are defined as sums over places $v$ of $K$ of local contributions $h_{\mathrm{Gr},K}^v(\cdot)$ and $h_{\mathrm{Tag},K}^v(\cdot)$ respectively.
Theorem \ref{itheo:main} is actually a consequence of parallelogram inequalities  for each of the local contributions  in $\hTag(\cdot)$.
Specifically, we prove the following result.
\begin{itheo}\label{itheo:local.tag}
Let $\phi:A\to K\{\tau\}$ be a Drinfeld module  of rank $r\geqslant 1$ defined over $K$, with stable reduction at all finite places of $K$. Let $v$ be a place of $K$. For any $G, H \in \KERS_K(\phi)$, we have
\begin{equation}\label{ieq:local.parineq}\tag{\text{\rotatebox{90}{$\lozenge$}}}
\hTag^v\big(\phi/(G\cap H)\big) + \hTag^v\big(\phi/(G+H)\big) \leqslant \hTag^v(\phi/G) + \hTag^v(\phi/H).
\end{equation}
\end{itheo}
Inequality~\eqref{ieq:local.parineq} is proved in Corollary \ref{prop:parineq.local.taguchi.finite} for finite $v$, and in Theorem~\ref{prop:parineq.local.taguchi.infinite} if $v$ is infinite.
If $v$ is infinite, it is worth noting  that inequality~\eqref{ieq:local.parineq} is actually an equality. The key technical result (Theorem \ref{prop:lattice_parallelogram_eq}) in this case, which we prove in section~\ref{sec:lattices}, is a relation between covolumes of $A$-lattices.
In case $v$ is finite, inequality~\eqref{ieq:local.parineq} is mainly based on a relation (Theorem \ref{thm:ineq.V.lcm.gcd}), proved in section~\ref{sec:valuation.polygons}, between valuation polygons of polynomials with coefficients in a valued field.

In the process of obtaining our main results, we also re-prove a general statement of independent interest, to the effect that isogenies preserve stable reduction at a place (Proposition \ref{prop:isog.stable.red}).
\medskip

We also obtain a parallelogram inequality for the local graded height at a finite place. 

\begin{itheo}\label{itheo:local.graded}
Let $\phi:A\to K\{\tau\}$ be a Drinfeld module  of rank $r\geqslant 1$ defined over $K$. Let $v$ be a \emph{finite} place of $K$.
For any $G, H \in \KERS_K(\phi)$, we have
\begin{equation}\label{ieq:local.parineq.graded}
\hG^v\big(\phi/(G\cap H)\big) + \hG^v\big(\phi/(G+H)\big) \leqslant \hG^v(\phi/G) + \hG^v(\phi/H).
\end{equation}
\end{itheo}

It is natural to wonder whether Theorem~\ref{itheo:local.graded} extends to the case of an infinite place $v$. We show that this is not the case.
Specifically, we discuss in section~\ref{sec:parineq.graded.infinite} a family of examples of parallelograms of Drinfeld modules of rank $2$. 
For different choices of parameters, we observe that the quantity
$\hG^v\big(\phi/(G\cap H)\big) + \hG^v\big(\phi/(G+H)\big) - \hG^v(\phi/G) - \hG^v(\phi/H)$ 
can be negative, zero, or positive.

\paragraph{General notation --}   
We fix a finite field $\F_q$ with cardinality $q$. 
We write $A\coloneq\F_q[T]$ for the ring of polynomials with coefficients in $\F_q$, and $K_0\coloneq\F_q(T)$ for the fraction field of $A$.
We let $K$ denote a finite field extension of $K_0$, \ie{}, a function field containing $\F_q$.
 
For any such field $K$, we let $M_{K}$ denote the set of places of $K$.
Each place $v\in M_K$ corresponds to a non-archimedean valuation $\ord_v : K \to \Z\cup\{\infty\}$.

The field $K_0$ has a distinguished place, denoted by $\infty$, which is characterized by the fact that $\ord_\infty(a) = -\deg(a)$ for any $a\in A$.
The other places of $K_0$ are called finite places: each finite place of~$K_0$ corresponds to a monic irreducible polynomial $P\in A$ and, for any $a\in A$, $\ord_v(a)$ is the multiplicity of $P$ as a divisor of $a$. 

A place $v\in M_K$ is called infinite if $v$ extends the place $\infty\in M_{K_0}$, and $v$ is called finite otherwise.
We let $M_K^\infty$ (resp. $M_K^f$) denote the set of infinite (resp. finite) places of $K$.
If a place $v\in M_K$ extends a place $w\in M_{K_0}$, we normalize $\ord_v$ so that $\ord_v(K^\ast) = \Z$. 
For any place $v\in M_K$, we write $K_v$ for the completion of $K$ at $v$.
To any place $v\in M_K$ lying above a place $w\in M_{K_0}$, we then associate its ramification index $e_v$, its residual degree $f_v$, and its local degree $n_v = [K_v : K_{0,w}]$. 
Recall that $n_v=e_v\,f_v$, see for instance Proposition 1.2.11 in \cite{BG06}.
Our choices of normalization are consistent with those in \cite{Papi, Gos12, DD99}, but differ from those in  \cite{BPR21}.

\numberwithin{equation}{section} 
\section{Isogenies between Drinfeld modules}

\subsection{Arithmetic of twisted polynomial rings}
Let $F$ be a field containing $\F_q$, and let $F^{\mathrm{sep}}$ denote the separable closure of $F$ in a chosen algebraic closure $\overline{F}$ of $F$.
We gather various useful facts about the  ring $F\{\tau\}$ of twisted polynomials in $\tau$: addition in $F\{\tau\}$  is the usual addition of polynomials and multiplication is given by
\[ \left(\sum_i f_i\,\tau^i\right)\cdot \left(\sum_j g_j\,\tau^j\right) = \sum_{k} \left(\sum_{i+j=k} f_i \,g_j^{q^i}\right)\tau^k.\]
The map that sends $\tau$ to $\mathrm{Frob}_q : x\mapsto x^q$ induces an $\F_q$-algebra isomorphism between $F\{\tau\}$ and the $\F_q$-algebra (endowed with composition) of endomorphisms of the additive group scheme $\mathbb{G}_{a, F}$ over $F$.
A polynomial $f\in F\{\tau\}$ will be called \emph{normalized} if its constant coefficient is $1$. 

For any $f=f_0\,\tau^0 + f_1\,\tau + \dots +f_d\,\tau^d \in K\{\tau\}$, we let
\[f(x) \coloneq f_0\,x + f_1\,x^q + \dots + f_d\, x^{q^d}\in F[x].\]
The map $x\mapsto f(x)$ is the $\F_q$-linear endomorphism of $\mathbb{G}_{a, F}$ associated to $f$.
The derivative (with respect to $x$) of $f(x)$ is the constant polynomial $f_0$.
When $f_0\neq 0$, the polynomial $f(x)$ is thus separable, \ie, has distinct roots in $\overline{F}$ (see \cite[p. 137]{Papi}).  
It is, in particular, the case if $f$ is normalized. 

\begin{definition}\label{defi:kernel}
For any separable polynomial $f\in F\{\tau\}$, we define the \emph{kernel of $f$} to be the set of roots in $\overline{F}$ of the corresponding polynomial $f(x)\in F[x]$:
\[ \ker(f) \coloneq \left\{ \beta\in\overline{F} \mid f(\beta)=0\right\} \subset\overline{F}. \]
\end{definition}
Since $x\mapsto f(x)$ is an $\F_q$-linear map, $\ker(f)$ is a finite dimensional $\F_q$-vector space. 
Note that $\card{\ker(f)} = \deg_x f(x) = q^{\deg_\tau f}$ and therefore $\dim_{\F_q}\ker(f) = \deg_\tau f$.
Conversely:
\begin{lem}\label{lem:ker.poly} 
Let $G$ be a finite dimensional $\F_q$-vector space contained in ${F}^{\mathrm{sep}}$, which is stable under the action of $\Gal(F^{\mathrm{sep}}/F)$.
Then there is a unique normalized  polynomial $f\in F\{\tau\}$ such that $\ker(f)=G$. 
\end{lem}
    \begin{proof}
    Given $G$ as in the statement and a non-zero constant $c\in F^\ast$, we let
    \[f_{G,c}(x) \coloneq c\, \prod_{\alpha\in G} (x-\alpha) \in F^{\mathrm{sep}}[x].\]
    Since $G\subset F^{\mathrm{sep}}$ is (globally) stable under the action of $\mathrm{Gal}(F^{\mathrm{sep}}/F)$, the polynomial $f_{G,c}(x)$ has coefficients in $F$. 
    The fact that $G$ is an $\F_q$-vector space allows to show that $f_{G, c}(x)\in F[x]$ is $\F_q$-linear, 
    \ie{}, that $f_{G,c}(x+\alpha\, y) = f_{G, c}(x) + \alpha\, f_{G,c}(y)$ for any $x,y\in\overline{F}$ and $\alpha\in\F_q$ (see the proof of \cite[Proposition 3.3.11]{Papi}).
    Hence $f_{G,c}(x)$  comes from a non-zero element $f_{G,c}\in F\{\tau\}$.
    
    It is clear that any separable polynomial $f\in F\{\tau\}$ with kernel $G$ is of the form $f_{G,c}$ for some $c\in F^\ast$. 
    There is exactly one choice of $c\in F^\ast$ such that $f_{G,c}$ has constant coefficient~$1$, namely $c= \prod_{\beta\in G\smallsetminus\{0\}}(-\beta)^{-1}$. 
    \end{proof}

Let $d,f\in F\{\tau\}$. One says that $d$ divides $f$ (on the right) if there exists $g\in F\{\tau\}$ such that $f= g\cdot d$. 
Equivalently, $d$ divides $f$ if and only if $f$ belongs to the left-ideal $F\{\tau\}\cdot d$.

A classical argument shows that the ring $F\{\tau\}$ admits a right-division algorithm (see \cite[Theorem 3.1.3]{Papi} or \cite[Proposition 1.6.2]{Gos12}, for instance), which immediately implies that any left-ideal in $F\{\tau\}$ is principal.

\begin{lem}\label{lemm:div.ker}
Let $d,f \in F\{\tau\}\smallsetminus\{0\}$ be separable polynomials. 
Then $d$ divides $f$ if and only if $\ker(d)\subset \ker(f)$.
\end{lem}
    \begin{proof} 
    It is clear that $\ker(d)\subset \ker(f)$ if $d$ divides $f$.
    Assume conversely that $\ker(d)\subset\ker(f)$.
    By the right-division algorithm in $F\{\tau\}$, there exists a unique pair $g, r\in F\{\tau\}$ such that $f = g\cdot d + r$ and $\deg_\tau r < \deg_\tau d$ or $r=0$.
    Since $r(\beta)= f(\beta) - g(d(\beta)) =0$ for all $\beta\in\ker(d)$, the polynomial $r(x)\in F[x]$ has at least $\card{\ker(d)}$ roots in $\overline{F}$.
    Since $d$ is separable, $\card{\ker(d)} = q^{\deg_\tau d}$.
    Having this many roots forces $r(x)$ to be~$0$.     
    \end{proof}

Given two separable polynomials $f, g\in F\{\tau\}$, it makes sense to define their (right)-GCD  as the unique normalized generator of the ideal $F\{\tau\}\cdot f + F\{\tau\}\cdot g $, and their (right)-LCM as the unique normalized generator of the ideal $F\{\tau\}\cdot f \cap F\{\tau\}\cdot g$.
Note that our terminology is non-standard, \cite[\S1.6]{Gos12} requires that the GCD be \emph{monic} instead.
The assumption that $f$ and $g$ be separable is necessary: for instance, $f=g=\tau$ do not admit a \emph{normalized} GCD.

\begin{lem}\label{lemm:ker.gcd.lcm}
    Let $f, g \in F\{\tau\}$ be non-zero separable polynomials. Write $d\in F\{\tau\}$ for their GCD and $\ell\in F\{\tau\}$ for their LCM. Then
    \[\ker(d) = \ker(f)\cap\ker(g) \qquad \text{and} \qquad \ker(\ell) = \ker(f)+\ker(g).\] 
\end{lem}
    \begin{proof}
    This follows from the definitions of GCD and LCM as well as the previous lemma.
    One can indeed rephrase the definition of the LCM as ``the LCM of $f$ and $g$ is the unique normalized polynomial $\ell$ of smallest degree such that $\ker(\ell)$ contains both $\ker(f)$ and $\ker(g)$''. A similar reformulation can be done for the GCD.
    \end{proof}

\subsection{Drinfeld modules and isogenies}\label{sec:Drinfeld_isogenies}
Let $F$ be a field extension of $K_0$. The field $F$ is then equipped with a (trivial) $A$-module structure {\em via} the inclusion $A\into K_0\into F$. 
In the terminology of \cite[\S3.2]{Papi}, $F$ is an $A$-field of $A$-characteristic $0$.

We recall the definition of a  Drinfeld module in this context:
\begin{definition} \label{def:D_module} 
  A Drinfeld module of rank $r \geqslant 1$ defined over $F$ is an $\F_q$-algebra homomorphism $\phi : A \to F\{\tau\}$, $a \mapsto \phi_a$ 
  such that, for each $a \in A$, 
  \begin{itemize}\setlength{\itemsep}{0em}\setlength{\topsep}{0em}      
      \item the constant coefficient of $\phi_a\in F\{\tau\}$ is $a\in A\subset F$, 
       \item and $\deg_\tau \phi_a = r \, \deg_T a$.
  \end{itemize}
\end{definition}

A Drinfeld module $\phi$ is entirely determined by the data $\phi_T\in F\{\tau\}$ because the $\F_q$-algebra ${A=\F_q[T]}$ is generated by $T$.
Such a morphism $\phi$ induces an $A$-module structure on $\overline{F}$ via the map $A\times\overline{F}\to \overline{F}$ given by $(a, x)\mapsto \phi_a(x)$. 
We refer to \cite[Chap.~12]{Rosen}, \cite[Chap.~4]{Gos12} or \cite[Chap.~3]{Papi}
for more details about the basics of the arithmetic theory of Drinfeld modules.

\begin{definition}\label{def:isog.drinfeld}
Given two Drinfeld modules $\phi, \psi$ defined over $F$, a morphism from $\phi$ to $\psi$ is a  polynomial $f\in F\{\tau\}$ such that 
\begin{equation}\label{eq:isogeny.definition}
     f\cdot\phi_a = \psi_a\cdot f \text{ for all } a\in A.
\end{equation}
Such a morphism will be denoted by $f:\phi\to\psi$.
A morphism is called an \emph{isogeny} if it is non-zero. 
We say that an isogeny is \emph{normalized} if its constant coefficient is $1$.
 
If there exists an isogeny between two Drinfeld modules, they will be called \emph{isogenous}.
\end{definition}

Since $A$ is generated by $T$ as an $\F_q$-algebra and since $\phi, \psi$ are $\F_q$-algebra morphisms, the ``commutation'' condition \eqref{eq:isogeny.definition} is equivalent to simply requiring that 
$f\cdot\phi_T = \psi_T\cdot f$. 

If there exists an isogeny $f:\phi\to\psi$, then $\phi$ and $\psi$ have the same rank (see \cite[Proposition~3.3.4(1)]{Papi}).
Moreover, in this case there exists a (non unique) \emph{dual isogeny} $\widehat{f}:\psi\to\phi$ (see \cite[Proposition~3.3.12]{Papi}). 
Therefore, ``being isogenous'' is an equivalence relation on the set of Drinfeld modules of rank~$r$ defined over $F$.
Note that an isogeny $f:\phi\to\psi$ is an isomorphism if and only if $\deg_\tau(f) = 0$, see \cite[\S3.8]{Papi}.
Actually, two Drinfeld modules $\phi, \psi$ of rank $r$ defined over $K$ are $K$-isomorphic if and only if there exists $c\in K^\ast$ such that $\psi = c^{-1}\cdot\phi\cdot c$.
    
\bigskip

Let $f:\phi\to\psi$ be an isogeny of Drinfeld modules over $F$.
Proposition 3.3.4(3) in \cite{Papi} implies that $f$ is automatically separable, \ie{}, $f(x)$ has distinct roots in $\overline{F}$.
We define, as we may, the kernel $\ker(f)$ of $f$ as the set of zeroes of $f(x)\in F[x]$ in some algebraic closure over $F$ (see Definition \ref{defi:kernel}):
\[\ker(f)=\left\{\beta\in\overline{F} \mid f(\beta)=0\right\}\subset \overline{F}.\]
Actually, since $f(x)$ is separable, we have $\ker(f)\subset F^{\mathrm{sep}}$.
Then $\ker(f)$ is a finite dimensional $\F_q$-vector space contained in $F^{\mathrm{sep}}$, with $\dim_{\F_q}\ker(f) = \deg_\tau(f)$.
Furthermore, notice that $\ker(f)$ is stable under the action of $\phi_a$ for all $a\in A$: indeed, for any $\beta\in\ker(f)$ and any $a\in A$, we have
\[f(\phi_a(\beta)) = (f\cdot\phi_a)(\beta) = (\psi_a\cdot f)(\beta) = \psi_a(f(\beta)) = \psi_a(0)= 0.\] 
Since $f\in F\{\tau\}$, the Galois action of $\Gal(F^{\mathrm{sep}}/F)$ on $F^{\mathrm{sep}}$ stabilizes $\ker(f)$.  
 
 This motivates the following definition:
\begin{definition} Let $\phi$ be a Drinfeld module over $F$.
    Let $\KERS_F(\phi)$ denote the set of finite dimensional $\F_q$-vector spaces $G$, contained in $F^{\mathrm{sep}}$, which are both stable under the action of $\mathrm{Gal}(F^\mathrm{sep}/F)$ and  stable under the action of~$\phi_a$ for all $a\in A$ (\ie{}, $\phi_a(x)\in G$ for all $x\in G$).
\end{definition}
Note that, for any $G_1, G_2 \in\KERS_F(\phi)$, both $G_1\cap G_2$ and $G_1+G_2$ are also elements of $\KERS_F(\phi)$.

\begin{prop}\label{prop:corres.ker.isog} 
    Let $\phi$ be a Drinfeld module defined over $F$.
\begin{enumerate}[(i)]\setlength{\itemsep}{0em}
    \item Let $f:\phi\to\psi$ be an isogeny from $\phi$ to a Drinfeld module $\psi$ defined over $F$. Then $\ker(f)\in\KERS_F(\phi)$.
    \item   For any $G\in \KERS_F(\phi)$, there is a unique Drinfeld module $\phi/G$ defined over $F$ and a unique \emph{normalized} isogeny $f_G:\phi\to\phi/G$ such that $\ker(f_G)=G$.
\end{enumerate}
\end{prop} 
For a given Drinfeld module $\phi$, the elements $G\in\KERS_F(\phi)$ therefore play a role similar to that of finite subgroup schemes of an abelian variety over $F$, in that they are both kernels of isogenies in their respective setting (see \cite{GLP26} for more details).
    \begin{proof} 
    The first item follows from the above discussion; we only need to prove the second item.
    Given an element $G\in \KERS_F(\phi)$ and a non-zero constant $c\in F^\ast$, we may define a polynomial $f_{G,c}\in F\{\tau\}$ as in the proof of Lemma \ref{lem:ker.poly}. 
    
    Consider the polynomial $\varpi \coloneq f_{G, c}\cdot\phi_T \in F\{\tau\}$. 
    For any $\alpha\in G$, we have $\phi_T(\alpha)\in G$ since we have assumed that $G$ is stable under the action of $\phi_T$. 
    Then $\varpi(\alpha) = f_{G, c}(\phi_T(\alpha)) = 0$. 
    In particular $G=\ker(f_{G, c})\subseteq \ker(\varpi)$.
    The link between kernels and divisibility in $F\{\tau\}$ (Lemma \ref{lemm:div.ker}) implies that there exists a unique $g_c\in F\{\tau\}$ such that $\varpi = g_c\cdot f_{G, c}$.
    The constant term of $f_{G, c}$ is non-zero, hence  the constant term of $g_c$ is $T$. We may now define a Drinfeld module $\psi_c : A\to F\{\tau\}$ by setting $\psi_{c} : T\mapsto g_c$. 
    By construction, we have $f_{G, c}\cdot \phi_T = \psi_{c,T}\cdot f_{G, c}$, \ie{}, $f_{G,c}\in F\{\tau\}$ is an isogeny $\phi\to\psi_c$.
    
    Notice that, for any $c\in F^\ast$, we have $f_{G, c} = c\, f_{G, 1}$ and  $\psi_c=c\,\psi_1\,c^{-1} \simeq \psi_1$. Therefore, up to isomorphim, there is a unique pair $(\psi, f)$ formed by a Drinfeld module $\psi$ defined over $F$ and an isogeny $f:\phi\to\psi$ with kernel $\ker(f)=G$. 
    Indeed any isogeny $f\in F\{\tau\}$ with kernel $G$ is of the form $f_{G,c}$ for some $c\in F^\ast$.
    Finally,  there is one choice of $c\in F^\ast$ such that $f_{G,c}$ has constant coefficient $1$, namely $c= \prod_{\beta\in G\smallsetminus\{0\}}(-\beta)^{-1}$. We set $\phi/G = \psi_c$ and $f_G=f_{G,c}$ for this choice of $c$.
    \end{proof}

We also need a factorization property of isogenies:
\begin{prop}\label{prop:facto.isog}
Let $\phi$ be a Drinfeld module defined over $F$. Assume there are isogenies $f_1:\phi\to\phi_1$ and $f_2:\phi\to\phi_2$ defined over $F$ such that $\ker(f_1)\subset\ker(f_2)$.
Then there exists a unique isogeny $g:\phi_1\to\phi_2$ such that $f_2=g\cdot f_1$. 
Moreover, if $f_1$ and $f_2$ are normalized, then so is $g$.
\end{prop}
    \begin{proof}
    By Lemma \ref{lemm:div.ker}, there exists a unique $g\in F\{\tau\}$ such that $f_2=g\cdot f_1$. 
    Using that $f_1:\phi\to\phi_1$ and $f_2:\phi\to\phi_2$ are isogenies we derive that, for any $a \in A$,
    \begin{align*}
        (g\cdot\phi_{1, a} - \phi_{2, a}\cdot g)\cdot f_1
        &= g\cdot\phi_{1, a}\cdot f_1 - \phi_{2, a}\cdot g\cdot f_1 
        = g\cdot f_1\cdot\phi_{1,a} - \phi_{2, a}\cdot f_2 \\
        &= f_2\cdot\phi_a - \phi_{2, a}\cdot f_2 
        =0.
    \end{align*}
    The ring $F\{\tau\}$ has no zero-divisors and $f_1\neq 0$, so we have $g\cdot\phi_{1, a} = \phi_{2, a}\cdot g$ for all $a\in A$, \ie{}, $g$ is an isogeny from $\phi_1$ to $\phi_2$.
    That $g$ is normalized if $f_1$ and $f_2$ are is immediate.
    \end{proof}

\subsection{Stable reduction of Drinfeld modules}\label{sec:stablered}
Let $K$ be a \emph{finite} field extension of $K_0$. Consider a Drinfeld module $\phi:A\to K\{\tau\}$ of rank $r\geqslant 1$ defined over $K$.
\begin{definition}
    For any  place $v$  of $K$, we say that \emph{$\phi$ has stable reduction at $v$} if there exists a Drinfeld module $\widetilde{\phi} : A\to K\{\tau\}$ defined over $K$ such that:
    \begin{itemize}\setlength{\itemsep}{0em}\setlength{\topsep}{0em}
        \item $\widetilde{\phi}$ is $K$-isomorphic to $\phi$, \ie{}, there exists $c\in K^\ast$ such that ${\widetilde{\phi}}_{\,T} = c^{-1}\cdot\phi_T\cdot c$,
        \item  and, writing ${\widetilde{\phi}}_{\,T}  = T + \sum_{i=1}^r g_i\, \tau^i$ with $g_i\in K$,  we have 
        $\min\{\ord_v(g_i), \, 1\leqslant i\leqslant r\} = 0$, \ie{}, the coefficients of ${\widetilde{\phi}}_{\,T}$ are integral at $v$ and at least one of them is a unit at $v$.
    \end{itemize}
    We say that $\phi$ has \emph{stable reduction everywhere} if it has stable reduction at all \emph{finite} places of $K$.
\end{definition}

See \cite[\S6.1]{Papi} or \cite[\S4.10]{Gos12} for  detailed treatments of this notion. 
We  note one important result about stable reduction (see \cite[Proposition 7.2]{Hayes} or \cite[Lemme 2.10]{DD99}):
\begin{theorem}
    Let $K$ be a finite extension of $K_0$, and let $\phi:A\to K\{\tau\}$ be a Drinfeld module of rank $r\geqslant 1$ defined over $K$.
    Then there exists a finite extension $K'/K$ such that the ``base changed'' Drinfeld module $\phi : A\to K'\{\tau\}$ has stable reduction everywhere.
\end{theorem}

Isogenies of Drinfeld modules preserve stable reduction at a place of $K$. This is the content of Proposition \ref{prop:isog.stable.red} below (see also Proposition 2.4 in \cite{Tag93}), which we defer to a later point in this paper because the proof of this statement requires tools developed in section \ref{sec:valuation.polygons}.

\section{\texorpdfstring{$A$}{A}-lattices and Drinfeld modules}\label{sec:lattices}

We let $\C_\infty$ denote the completion of an algebraic closure of the completion $K_{0, \infty}$ of $K_0$ at the place~$\infty$.
The field $\C_\infty$ is equipped with a non-archimedean valuation, extending $\ord_\infty:K_0^\ast\to\Z$. We let $|\cdot|_\infty$ denote the corresponding normalized absolute value, and we endow $\C_\infty$ with the induced topology.

\subsection{\texorpdfstring{$A$}{A}-lattices and their isogenies}
We recall the following definition (see \cite[Def. 4.1.2]{Gos12}  or \cite[Def. 5.2.1]{Papi}).

\begin{definition} \label{def:lattice}
An \emph{$A$-lattice} (or simply a \emph{lattice}) is a finitely generated $A$-submodule $\Lambda$ of $\C_\infty$, 
which is discrete as a subset of $\C_\infty$ (\ie{}, $\Lambda$ has finite intersection with any ball in $\C_\infty$).

An $A$-lattice $\Lambda$ is necessarily torsion-free; the \emph{rank} of $\Lambda$ is its rank as a finitely generated $A$-module.
\end{definition}

Let $\Lambda$ be an $A$-lattice of rank $r\geqslant 1$. Since $A=\F_q[T]$ is a principal ideal domain and since $\Lambda\subset\C_\infty$ is discrete, there exist $\omega_1, \omega_2, \dotsc, \omega_r \in \C_\infty$ such that $\Lambda = A\,\omega_1 +A\,\omega_2 +\dotsb +A\,\omega_{r}$ and the $\omega_i$ are $K_{0,\infty}$-linearly independent (see \cite[Proposition 5.3.8]{Papi}).

\begin{definition} Let $\Lambda_1, \Lambda_2\subset \C_\infty$ be two $A$-lattices \emph{of the same rank} $r\geqslant 1$.
A morphism $c:\Lambda_1\to\Lambda_2$  is an element $c\in\C_\infty$ such that $c\,\Lambda_1\subset \Lambda_2$. Such a morphism is called an \emph{isogeny} if $c\neq 0$.

If there exists an isogeny between two lattices, they will be called \emph{isogenous}.    
\end{definition}\label{def:lattice.morphism}
Two lattices $\Lambda_1$ and $\Lambda_2$ are isomorphic if and only if there is a $c \in \C_\infty\smallsetminus\{0\}$ such that $c\,\Lambda_1 = \Lambda_2$.

For any isogeny $c:\Lambda_1\to\Lambda_2$, the $A$-lattice $c\, \Lambda_1$ has finite index in $\Lambda_2$, \ie{}, the quotient $A$-module $\Lambda_2/c\,\Lambda_1$ is a finite $A$-module. 
Since $A=\F_q[T]$ is a principal ideal domain, the structure theorem for torsion $A$-modules (see \cite[III.7, Thm. 7.7]{Lang}) implies that $\Lambda_2/c\,\Lambda_1$ is a  direct sum of finitely many $A$-modules of the form $A/a_i\,A$ (with non-constant $a_i\in A$).
In particular, the quotient $\Lambda_2/c\,\Lambda_1$ is a finite dimensional $\F_q$-vector space and thus a finite set: we usually denote its cardinality $|\Lambda_2/c\,\Lambda_1|$ (a power of $q$) as an index $\big(\Lambda_2:c\,\Lambda_1\big)$.

 We refer to Theorem 5.3.10 of \cite{Papi} for many more details about $A$-lattices and isogenies of $A$-lattices.

\subsection{Covolume of \texorpdfstring{$A$}{A}-lattices}
Let $\Lambda$ be an $A$-lattice of rank $r\geqslant 1$. 
A \emph{successive minimum basis} of $\Lambda$ is an ordered basis $(\eta_1, \dots, \eta_r)$ of $\Lambda$ as an $A$-module such that $|\eta_1|_\infty\geqslant |\eta_2|_\infty\geqslant \dots \geqslant |\eta_r|_\infty,$ 
and, for any $k=1, \dots, r$, 
\[ |\eta_k|_\infty = \min\left(\left\{ \left|\lambda - {\textstyle\sum_{i=k+1}^r a_i\, \eta_i} \right|_\infty , \ (a_{k+1}, \dots, a_r)\in A^{r-k+1},\,  \lambda\in\Lambda \right\}\smallsetminus\{0\}\right). \]
In other words, $\eta_r$ has minimal absolute value among non-zero elements of $\Lambda$, and each $\eta_k$ has minimal absolute value among elements of 
$\Lambda \smallsetminus \big(A\,\eta_{k+1} + \dots + A\,\eta_r\big)$. 
See \cite[\S4.1]{BPR21} or \cite[\S5.5]{Papi}.

By the discussion in \cite[\S5.5]{Papi}, any $A$-lattice $\Lambda$ of positive rank admits a successive minimum basis; 
such a basis is not necessarily unique, but the tuple $(|\eta_1|_\infty, \dots, |\eta_r|_\infty)$ is an invariant of $\Lambda$.
We then define the \emph{covolume of $\Lambda$} by 
\begin{equation}
    \covol(\Lambda) \coloneq |\eta_1|_\infty\,|\eta_2|_\infty\,\cdots\,|\eta_r|_\infty,  
\end{equation}
for any successive minimum basis $(\eta_1, \dots, \eta_r)$ of $\Lambda$. 
By the previous paragraph, $\covol(\Lambda)$ is well-defined, \ie{}, independent of the choice of a successive minimum basis of $\Lambda$.

It will be relevant for us that the covolume interacts well with isogenies of $A$-lattices: 
\begin{lem}\label{lem:covol.isog}
Let $c\in\C_\infty^\ast$ be an isogeny $c:\Lambda\to\Lambda'$ between two $A$-lattices $\Lambda, \Lambda'$ of rank $r\geqslant 1$.
Then we have
\[ \covol(\Lambda')  =  |c|_\infty^{r}\,\big(\Lambda':c\,\Lambda\big) \,  \covol(\Lambda). \]
In particular, if $\Lambda$ and $\Lambda'$ are two $A$-lattices of rank $r\geqslant 1$ such that $\Lambda \subseteq\Lambda'$, then
\begin{equation*}
    \covol(\Lambda')= \big(\Lambda':\Lambda \big) \,\covol(\Lambda).
\end{equation*}
\end{lem}
    \begin{proof}
        This follows directly from Lemma 4.1 in \cite{BPR21}.
    \end{proof}

\subsection{Drinfeld modules and lattices}\label{sec:uniformization}
We now recall the link between the theory of $A$-lattices and that of Drinfeld modules over $\C_\infty$.
\begin{theorem}\label{theo:corresp.drinfeld.lattices}
    For any $r\geqslant 1$, there is an equivalence of categories between
    \begin{enumerate}[(a)]\setlength{\itemsep}{0em}
        \item the category whose objects are the Drinfeld modules $\phi : A\to\C_\infty\{\tau\}$ of rank $r$ defined over $\C_\infty$, and whose morphisms are the morphisms of Drinfeld modules as in Definition \ref{def:isog.drinfeld},
        \item the category whose objects are $A$-lattices $\Lambda\subset \C_\infty$ of rank $r$, and whose morphisms are the morphisms of lattices as in Definition \ref{def:lattice.morphism}.
    \end{enumerate}
\end{theorem}

This is Theorem 4.6.9 in \cite{Gos12}. The correspondence can be made more explicit by using the exponential function, as we now briefly explain.
To any $A$-lattice $\Lambda\subset\C_\infty$, one associates a surjective $\Lambda$-periodic entire function $\ee_\Lambda:\C_\infty\to\C_\infty$ defined by a convergent product
\[\forall z\in\C_\infty, \qquad \ee_\Lambda(z) = z\, \prodp_{\lambda \in \Lambda} \left(1-\frac{z}{\lambda}\right), \]
where $\prod'$ means that the product is taken over all non-zero elements of $\Lambda$. We refer to \cite[\S4.2]{Gos12} or \cite[\S5.1]{Papi} for details.
Let $\phi^\Lambda_0(x)\coloneq 0$ and, for any $a\in A\smallsetminus\{0\}$, let
\[ \phi^\Lambda_a(x) \coloneq a^{-1}\,\ee_\Lambda\big( a\, \ee_\Lambda^{-1}(x)\big), \]
where $\ee_\Lambda^{-1}(x)$ denotes an element $z\in\C_\infty$ such that $\ee_\Lambda(z)=x$. 
One shows that $\phi^\Lambda_a(x)$ is actually an $\F_q$-linear polynomial in $\C[x]$, hence comes from an element $\phi^\Lambda_a\in \C_\infty\{\tau\}$.
One then proves that the map, $\phi^\Lambda : A\to\C_\infty\{\tau\}$ given by $a\mapsto\phi^\Lambda_a$ 
defines a Drinfeld module of rank $r$. 
The next step is to show that any Drinfeld module of rank $r$ over $\C_\infty$ arises in this way.
We refer to  \cite[\S4.3]{Gos12} or \cite[\S5.2]{Papi} for more details and proofs. 

On the level of morphisms, the correspondence is given as follows. Let $\phi, \psi$ be two Drinfeld modules defined over $\C_\infty$. 
Write $\Lambda_\phi, \Lambda_\psi$ for the $A$-lattices associated to $\phi$ and $\psi$ by the previous paragraph.
Then the map
\[\begin{array}{ccc}
\{\text{morphisms of Drinfeld modules } \phi\to \psi\} &\longrightarrow& \{\text{morphisms of $A$-lattices }\Lambda_\phi\to\Lambda_\psi\}  \\
     f\in\C_\infty\{\tau\} &\longmapsto& \text{constant coefficient of } f
\end{array} \]
is the desired bijection.
The reader is invited to consult \cite[\S4.6]{Gos12} or \cite[\S5.2]{Papi} for more details.
For later use, we record the following result, which is essentially Theorem 5.2.10 in \cite{Papi}:
\begin{prop}\label{prop:corresp.isog.lattices}
Let $f:\phi\to\psi$ be an isogeny of Drinfeld modules defined over $\C_\infty$. Write $f_0\in\C_\infty$ for the constant coefficient of $f$.
Let $\Lambda_\phi$ and $\Lambda_\psi$ be the $A$-lattices associated to $\phi$ and $\psi$ by Theorem~\ref{theo:corresp.drinfeld.lattices}, respectively.
The isogeny $c :\Lambda_\phi \to\Lambda_\psi$ of lattices induced by $f$ is $c=f_0$.
Moreover, there is an isomorphism of $\F_q$-vector spaces $\ker(f) \simeq \Lambda_\psi / f_0\, \Lambda_\phi$ and we have $\deg_\tau f = \log_q\big(\Lambda_\psi : f_0\, \Lambda_\phi\big)$.
\end{prop}
In particular,  a \emph{normalized} isogeny $f:\phi\to\psi$ induces an \emph{inclusion} of $A$-lattices $\Lambda_\phi\subseteq\Lambda_\psi$, and 
\[  \deg_\tau f = \dim_{\F_q}\ker(f) = \dim_{\F_q} \Lambda_\psi/\Lambda_\phi = \log_q\big(\Lambda_\psi :  \Lambda_\phi\big).\]

\subsection{A covolume identity}
We now prove the following:
\begin{theorem}\label{prop:lattice_parallelogram_eq}
Let $\Lambda, \Lambda_1, \Lambda_2$ be three $A$-lattices of the same rank $r\geqslant 1$ such that $\Lambda \subseteq \Lambda_1\cap\Lambda_2$. 
Then $\Lambda_1\cap\Lambda_2$ and $\Lambda_1+\Lambda_2$ are also $A$-lattices of rank $r$.
Moreover, we have
    \begin{equation}\label{eq:lattice.parallelogram.ineq}
    \log_q \covol(\Lambda_1\cap\Lambda_2)  + \log_q \covol(\Lambda_1 + \Lambda_2) = \log_q \covol(\Lambda_1) + \log_q \covol(\Lambda_2). 
    \end{equation}
\end{theorem}

Identity \eqref{eq:lattice.parallelogram.ineq} may be thought of as a ``parallelogram equality'' for $A$-lattices.
The data of the statement indeed naturally gives rise to a diagram of $A$-lattices where arrows represent inclusions:
\begin{center} 
\begin{tikzcd}
                           & \Lambda_1\cap\Lambda_2 \arrow[ld, hook] \arrow[rd, hook] &                            \\
\Lambda_1 \arrow[rd, hook] &                                                          & \Lambda_2\,. \arrow[ld, hook] \\
                           & \Lambda_1+\Lambda_2                                      &                           
\end{tikzcd}
\end{center}
    \begin{proof}
    Recall that $A=\F_q[T]$ is a principal ideal domain.
    As a sub-$A$-module of the finitely generated free $A$-module $\Lambda_1$ of rank $r$, $\Lambda_1\cap\Lambda_2$ is also a free $A$-module of rank $\leqslant r$.
    On the other hand, the $A$-module $\Lambda$ of rank $r$ is contained in $\Lambda_1\cap\Lambda_2$. The latter must therefore have rank exactly $r$.  
    In particular,  $\Lambda_1\cap\Lambda_2$  has finite index in $\Lambda_2$. 
    The sum $\Lambda_1+\Lambda_2$ is a sub-$A$-module of the free $A$-module $\C_\infty$, and it contains the $A$-module $\Lambda_1$ of rank $r$. 
    Hence $\Lambda_1+\Lambda_2$ is a free finitely generated $A$-module of rank $\geqslant r$.
    The ``second isomorphism theorem'' for $A$-modules (see p. 120 in \cite[\S{}III.1]{Lang}) shows that there is an isomorphism of $A$-modules
    \begin{equation}\label{eq:third.isom.thm}
        (\Lambda_1+\Lambda_2)/\Lambda_1 \simeq \Lambda_2/(\Lambda_1\cap\Lambda_2).
    \end{equation}
    In particular, the quotient $(\Lambda_1+\Lambda_2)/\Lambda_1$ is a finite $A$-module, and $\Lambda_1+\Lambda_2$ has rank equal to $r$.
    The fact that both $\Lambda_1\cap \Lambda_2$ and $\Lambda_1+\Lambda_2$ are discrete subsets of $\C_\infty$ follows easily from the discreteness of $\Lambda_1$, $\Lambda_2$ in $\C_\infty$.
    The first assertion is proved.
    
    Combined with \eqref{eq:third.isom.thm}, the ``third isomorphism theorem'' for $A$-modules (\textit{loc. cit.}) yields $A$-module isomorphisms
    \[\frac{(\Lambda_1 +\Lambda_2)/\Lambda}{\Lambda_1/\Lambda} \simeq \frac{\Lambda_1+\Lambda_2}{\Lambda_1} \simeq \frac{\Lambda_2}{\Lambda_1\cap \Lambda_2} \simeq \frac{\Lambda_2/\Lambda}{(\Lambda_1\cap \Lambda_2)/\Lambda}.\]
    In particular, we obtain
    $\big(\Lambda_1+\Lambda_2 :\Lambda\big) \, \big(\Lambda_1\cap\Lambda_2 : \Lambda\big) = \big(\Lambda_1:\Lambda\big)\,\big(\Lambda_2:\Lambda\big)$.
    Multiplying on both sides by~$\covol(\Lambda)^2$ and using the second identity in Lemma \ref{lem:covol.isog}, we get
    \[\covol(\Lambda_1+\Lambda_2)\,\covol(\Lambda_1\cap\Lambda_2) = \covol(\Lambda_1)\,\covol(\Lambda_2).\]
    Taking the image by the logarithm $\log_q$ yields the desired equality.
    \end{proof}

\section{Valuation polygons of twisted polynomials}\label{sec:valuation.polygons}
For the duration of this section, we work in the following setting:
let $F$ be a field containing $\F_q$, equipped with a discrete non-archimedean valuation $\nu:F\to\Z\cup\{\infty\}$, normalized so that $\nu(F^\ast) = \Z$.
We fix an algebraic closure $\overline{F}$ of $F$ and choose an extension of $\nu$ to $\overline{F}$, also denoted by $\nu:\overline{F}\to\Q\cup\{\infty\}$.

Definitions and basic properties of the objects described in this section can essentially be found in \cite[\S1]{Tag91} and \cite[\S2]{Tag93}. 
We thought it worthwhile to rewrite detailed definitions and some of the proofs in our notation, in order to keep this paper self-contained and clarify some steps.

\subsection{Valuation polygons}\label{sec:val.polygon}

Let $\displaystyle{f=\sum_{i=0}^d f_i\,\tau^i \in F\{\tau\}}$ be a separable polynomial.
For $z\in\R$, let
\begin{equation}\label{eq:def.val.pol}
    \V^\nu_f(z) \coloneq \min\{ \nu(f_i) + q^i\, z, \   0\leqslant i \leqslant d\}.
\end{equation}
The map thus defined is continuous and (strictly) increasing on $\R$, and $\V^\nu_f:\R\to\R$ is a bijection.
By definition, its graph is a concave broken line with finitely many slopes, each of which is an element of $\{1, q, \dotsc, q^d\}$.
A direct computation shows that the abscissae of the break-points are exactly the elements of 
\[\mathcal{E}_f \coloneq \left\{ - \frac{\nu(f_j)-\nu(f_i)}{q^j - q^i},  \ 0\leqslant i\neq j\leqslant d\right\}.\] 
One checks that 
$\{\nu(\beta), \, \beta\in\ker(f)\smallsetminus\{0\}\} \subseteq \mathcal{E}_f \subset \Q$.

We also note another way of characterizing $\V^\nu_f$:
\begin{lem}\label{lem:vee.f.ord}
    For all $x\in\overline{F}$, we have $\nu\big(f(x)\big) \geqslant \V^\nu_f\big(\nu(x)\big)$. Furthermore, for all $x\in\overline{F}$ such that $\nu(x)\notin \mathcal{E}_f$, we have $\nu\big(f(x)\big) = \V^\nu_f\big(\nu(x)\big)$.
\end{lem}
    \begin{proof}
        The first statement is a direct consequence of the non-archimedean triangle inequality (see \cite[\S2.1]{Papi} for instance).
        By definition of the set $\mathcal{E}_f$, for any $x\in\overline{F}$ with $\nu(x) \notin \mathcal{E}_f$, the $d+1$ numbers $\nu(f_i) + q^i\, \nu(x) = \nu(f_i\, x^{q^i})$ ($0\leqslant i\leqslant d$) are pairwise distinct. 
        The second assertion now follows from the first one, and the case of equality in the triangle inequality.
    \end{proof}

The graph of $\V^\nu_f$ is sometimes called the \emph{valuation polygon of $f$}, and it is in a sense ``dual'' to the Newton polygon of $f(x)/x\in F[x]$, as explained in \cite[\S2.1--\S2.2]{Robba}.
Write the elements of $\mathcal{E}_f$ in increasing order: $e_1< e_2 < \dots < e_m$. 
It will be convenient for later to rename $e_m$: 
\begin{equation}\label{eq:def.lambda}
\lambda_\nu(f) \coloneq  -\min\left\{\frac{\nu(f_i)-\nu(f_0)}{q^i-1} \,,\ 1\leqslant i\leqslant d\right\}
= \max\{\nu(\beta), \, \beta\in\ker(f)\smallsetminus\{0\}\}.
\end{equation}
It is clear  that $\V^\nu_f(z)\leqslant z +\nu(f_0)$ for all $z\in\R$. Notice (as in \cite[Remark 1.4]{Tag91}) that
\begin{equation}
    \lambda_\nu(f)=\min\left\{ z\in \R \,\mid\, \V^\nu_f(z) = z + \nu(f_0)\right\}.
\end{equation}

For any $z\in\R$, set
\[\N^\nu_f(z) \coloneq \card{\{ \beta\in \ker(f) \mid \nu(\beta)\geqslant z\}}.\]
The map $z\mapsto \N^\nu_f(z)$ is non-increasing, left-continuous and piecewise constant.
Since $f(x)$ is a \linebreak {$\F_q$-linear} polynomial and $\nu$ is non-archimedean, the set $\left\{ \beta\in \ker(f) \mid \nu(\beta)\geqslant z\right\}$ is a finite dimensional $\F_q$-vector space, so that $\N^\nu_f$ has values in $\{1, q, \dotsc, q^d\}$. 
By definition, the graph of $\N^\nu_f$ has jumps at elements of~$\mathcal{E}_f$ and is constant on each segment $(e_k, e_{k+1}]$, as well as on $(-\infty, e_1]$ and on $(e_m, +\infty)$.
For $z \leqslant e_1$, it is clear that $\N^\nu_f(z) = \card{\ker(f)} = q^d$.
For $z > e_m$, we have $\N^\nu_f(z) = 1$ (the only root of $f(x)$ with valuation larger than $e_m$ is $0$).
For $1 \leqslant k < m$ and $z\in(e_k, e_{k+1}]$, one readily checks (by an argument similar to the proof of  Theorem 2.5.2 in \cite{Papi} for instance, see \cite[\S2.1]{Robba}) that $\V^\nu_f(z) = \nu(f_i) + \N^\nu_f(z)\,z$, \ie{}, that 
the index $i\in\{1, \dots, r\}$ realizing the minimum in $\V^\nu_f(z)$ is $\log_q \N^\nu_f(z)$.
We conclude that the maps $\V^\nu_f$ and $\N^\nu_f$ are related as follows:
\[ \V^\nu_f(z) = \begin{cases}
    \nu(f_d) + q^d\, z  & \text{if }z < e_1, \\
    \nu(f_i)  + q^i \,z &\text{if } z\in(e_k, e_{k+1}] \text{ for }1 \leqslant k < m, \ \text{where} \ i = \log_q \N^\nu_f(z), \\
    \nu(f_0) + q^0\, z  & \text{if }z > e_m = \lambda_\nu(f).
\end{cases} \]
In other words, we have 
\begin{lem}
For any $z\in\R\smallsetminus \mathcal E_f$, the map $\V^\nu_f$ is differentiable at $z$, and 
$(\V^\nu_f)'(z) = \N^\nu_f(z)$.
\end{lem}

\subsection{An inequality}
We now  prove the main result of this section, which relates the valuation polygons of two normalized polynomials to those of their GCD and LCM.

\begin{theorem} \label{thm:ineq.V.lcm.gcd}
Let $f$ and $g\in F\{\tau\}$ be normalized polynomials. Let $d$ and $\ell$ denote their normalized GCD and LCM in $F\{\tau\}$, respectively. 
Then, for all $z\in\R$, 
\[ \V^\nu_d(z) +\V^\nu_\ell(z) \leqslant \V^\nu_f(z) +\V^\nu_g(z). \]
\end{theorem}
In the special case $d=1$ (\ie{}, when $f$ and $g$ have trivial GCD in $F\{\tau\}$) this inequality becomes $\V^\nu_\ell(z) \leqslant \V^\nu_f(z) +\V^\nu_g(z) -z$ for all $z\in\R$.
Notice indeed that $\V^\nu_1(z)=z$ for all $z\in\R$. 

We begin by proving a lemma: 
\begin{lem} \label{lem:Nf.parallelogram}
    In the context of the theorem, we have  
    \[\forall z\in\R, \qquad   \N^\nu_d(z)\,\N^\nu_\ell(z) \geqslant \N^\nu_f(z)\,\N^\nu_g(z) \quad\text{ and }\quad \N^\nu_d(z) +\N^\nu_\ell(z) \geqslant \N^\nu_f(z) +\N^\nu_g(z). \]
\end{lem}
    \begin{proof} 
    For  $h\in \{f,g,d,\ell\}$, we write $S_h(z) \coloneq \{\beta\in\ker(h) \mid \nu(\beta)\geqslant z\}$.
    Notice that $S_h(z)$ is an $\F_q$-vector space contained in $\ker(h)$, and that $\card{S_h(z)} = \N^\nu_h(z)$.
    Recall from Lemma \ref{lemm:ker.gcd.lcm} that $\ker(\ell) = \ker(f) + \ker(g)$ and $\ker(d) = \ker(f) \cap \ker(g)$.
    For a fixed $z \in \R$, consider the map 
    \[ \sigma: S_f(z)\times S_g(z) \to S_\ell(z), \quad (\beta, \beta') \mapsto \beta+\beta'. \]
    This defines an $\F_q$-linear map because $\ker(f)+\ker(g) = \ker(\ell)$ and, 
    for any $\beta, \beta'$ such that $\nu(\beta)\geqslant z$, $\nu(\beta')\geqslant z$, we have
    $\nu(\beta+\beta')\geqslant\min\{\nu(\beta),\, \nu(\beta')\}\geqslant z$.
    The kernel of $\sigma$ is a diagonally embedded copy of $S_f(z)\cap S_g(z) = S_d(z)$ in $S_f(z)\times S_g(z)$.
    Hence $\sigma$ induces an injective $\F_q$-linear map \[\big(S_f(z)\times S_g(z)\big) / S_d(z) \longrightarrow S_\ell(z).\] 
    Comparing cardinalities of the domain and codomain of this latter map yields the first inequality.
    
    To prove the second inequality, note  that both $\N^\nu_f(z)$ and $\N^\nu_g(z)$ are contained in the interval $\big[\N^\nu_d(z),  \N^\nu_\ell(z) \big]$ for all $z\in\R$; thus
    \begin{align*}
        \N^\nu_d(z)\,\N^\nu_\ell(z) &\geqslant \N^\nu_f(z)\,\N^\nu_g(z) \\
        &= (\N^\nu_f(z)-\N^\nu_d(z)) (\N^\nu_g(z)-\N^\nu_d(z)) +\N^\nu_f(z)\,\N^\nu_d(z) +\N^\nu_g(z)\,\N^\nu_d(z) -  \N^\nu_d(z)^2 \\
        &\geqslant \N^\nu_f(z)\,\N^\nu_d(z) + \N^\nu_g(z)\,\N^\nu_d(z) - \N^\nu_d(z)^2.
    \end{align*}
    The result follows by dividing by $\N^\nu_d(z)\geqslant 1$.
    \end{proof}

    \begin{proof}[Proof (of Theorem \ref{thm:ineq.V.lcm.gcd})] 
    Given  $z\in\R$, we pick $z_0 = 1+ \max\{ z, \lambda_\nu(f), \lambda_\nu(g), \lambda_\nu(d),\lambda_\nu(\ell)\}$. 
    
    For  $h\in \{f,g,d,\ell\}$,
    recall that $\V^\nu_h$ is continuous on~$\R$, is differentiable almost everywhere (the exceptional abscissae are the elements of $\mathcal{E}_h$),  and that $(\V^\nu_h)'(t)=\N^\nu_h(t)$ for any $t\in\R\smallsetminus\mathcal{E}_h$. Also note that $\V^\nu_h(z_0)=z_0$ because $h$ is normalized.
    We thus have
    \[ \V^\nu_h(z) = z_0 - \int_z^{z_0} \N^\nu_h(t)\dd t. \]
    Using that $z < z_0$, we now deduce from the above lemma that
    \begin{align*}
        \V^\nu_d(z) + \V^\nu_\ell(z) 
        &= 2\,z_0 - \int_z^{z_0}\left( \N^\nu_d(t)+\N^\nu_\ell(t)\right)\dd t 
        \leqslant 2\,z_0 - \int_z^{z_0} \left( \N^\nu_f(t) + \N^\nu_g(t) \right)\dd t \\
        &= 2\,z_0 - \V^\nu_f(z_0)+\V^\nu_f(z) - \V^\nu_g(z_0)+\V^\nu_g(z)
        = \V^\nu_f(z) + \V^\nu_g(z).
    \end{align*}
    The desired inequality is proved.
    \end{proof}

\subsection{Another property of valuation polygons}
In order to use the results of this section in the context of Drinfeld modules (see \S\ref{sec:isog.local.height} below), we require  Proposition \ref{prop:lambda.poly.isog} below. 

We begin by a lemma (which is essentially \cite[Lemma 1.2]{Tag91}):  

\begin{lem}\label{lemm:V.compatibility.composition}
Let $f,g\in F\{\tau\}$ be separable polynomials.
For all $z\in\R$, we have
\[\V^\nu_{f\cdot g}(z) = \V^\nu_f\circ \V^\nu_g(z).\]    
\end{lem}
    \begin{proof}
    By Lemma \ref{lem:vee.f.ord}, for almost all $z\in\Q$ and any $x\in\overline{F}$ with $\nu(x) = z$, we have
    \[ \V^\nu_{f\cdot g}(z) = \V^\nu_{f\cdot g}(\nu(x)) 
        = \nu\big((f\cdot g)(x)\big)
        = \nu\big(f(g(x))\big)
        = \V^\nu_f(\nu(g(x)) )
        =\V^\nu_f(\V^\nu_g(z)).\]
    The two maps $z\mapsto \V^\nu_{f\cdot g}(z)$ and $z\mapsto \V^\nu_f\circ \V^\nu_g(z)$ are continuous on $\R$. We have just noted that they coincide on a dense subset of $\R$, so they are equal.
    \end{proof}

The main result of this paragraph is the following statement, generalizing \cite[Lemma 2.3]{Tag93}.
\begin{prop} \label{prop:lambda.poly.isog}
    Let $\phi$, $\psi$, $f \in F\{\tau\}$ be separable polynomials such that $f \cdot \phi = \psi \cdot f$. Then we have  $\phi_0 = \psi_0$, and
    \begin{itemize}
        \item If $\nu(\psi_0) = \nu(\phi_0) \geqslant 0$, then $\lambda_\nu(\psi) \leqslant \V^\nu_f(\lambda_\nu(\phi))$.
        \item If $\nu(\psi_0) = \nu(\phi_0) \leqslant 0$, then $\lambda_\nu(\psi) \geqslant \V^\nu_f(\lambda_\nu(\phi))$.
    \end{itemize}
    In particular, if $\nu(\phi_0) = \nu(\psi_0) = 0$, then $\lambda_\nu(\psi) = \V^\nu_f(\lambda_\nu(\phi))$.
\end{prop}
    \begin{proof}
    The fact that $\psi_0 = \phi_0$ follows from comparing the coefficient of $\tau^0$ on both sides of the relation $f \cdot \phi = \psi \cdot f$, since $f_0 \neq 0$.
    Recall that, for a separable polynomial $h\in F\{\tau\}$, $\V^\nu_h$ is a strictly increasing bijection and that, for any $z\in\R$, we have $z \geqslant \lambda_{\nu}(h) \iff \V^\nu_h(z) = z +\nu(h_0)$. 
    Given the equality $f \cdot \phi = \psi \cdot f$, the above lemma yields $(\V^\nu_f)^{-1}\circ\V^\nu_\psi = \V^\nu_\phi\circ (\V^\nu_f)^{-1}$.
    In this proof, we write $\mu \coloneq \V^\nu_f(\lambda_\nu(\phi))$.  
    
    We first assume that $\nu(\psi_0) = \nu(\phi_0) \geqslant 0$. For any $z\in\R$, we let $y = (\V^\nu_f)^{-1}(z)$. We have 
    \begin{equation*}
    z\geqslant \mu 
    \iff y \geqslant \lambda_\nu(\phi)
    \iff \V^\nu_\phi (y) = y + \nu(\phi_0)
    \iff \V^\nu_\psi (z) = \V^\nu_f(y + \nu(\phi_0)).
    \end{equation*}
    If $z\geqslant \mu$, we  thus have 
    $\V^\nu_\psi(z) = \V^\nu_f(y+\nu(\phi_0)) \geqslant \V^\nu_f(y) +\nu(\phi_0) = z +\nu(\psi_0)$,
    since $(\V^\nu_f)'(t)\geqslant 1$ for almost all $t\in\R$.
    So for any $z\geqslant\mu$, we have $\V^\nu_\psi(z)\geqslant z+\nu(\psi_0)$. This implies that $z\geqslant \lambda_\nu(\psi)$. 

    Next we assume that $\nu(\psi_0) = \nu(\phi_0) \leqslant 0$. For any $z\in\R$, we let $y = (\V^\nu_f)^{-1}(z)$. We have
   \begin{equation*}
   z \geqslant \lambda_\nu(\psi) 
   \iff \V^\nu_\psi(z) = z + \nu(\psi_0) 
   \iff \V^\nu_f \circ \V^\nu_\phi(y) = z + \nu(\psi_0)
   \iff \V^\nu_\phi(y) = (\V^\nu_f)^{-1}(z+ \nu(\psi_0)). 
    \end{equation*}
    Hence, if $z\geqslant \lambda_\nu(\psi)$, then 
    \[\V^\nu_\phi(y) = (\V^\nu_f)^{-1}(z+ \nu(\psi_0)) \geqslant (\V^\nu_f)^{-1}(\V^\nu_f(y)) + \nu(\psi_0) = y + \nu(\phi_0),\]
    because the derivative of $(\V^\nu_f)^{-1}$ is bounded from above by $1$.  
    Therefore  we have $\V^\nu_\phi(y) \geqslant y + \nu(\phi_0)$ for any $z\geqslant\lambda_\nu(\psi)$. This implies that $y = (\V^\nu_f)^{-1}(z)\geqslant \lambda_\nu(\phi)$, which means that $z\geqslant \mu$. 
    \end{proof}

\begin{ex}
If $\nu(\phi_0)=\nu(\psi_0)$ is non zero in the above proposition, it is possible that the inequalities between $\lambda_\nu(\psi)$ and $\V^\nu_f(\lambda_\nu(\phi))$ are strict. We give an example.

Let $q\geqslant 2$ be a prime power and let $F=\F_q(T)$ equipped with the discrete valuation $\nu = -\deg$ (corresponding to the place at infinity). 
For any $a\in A\smallsetminus\{0\}$, let 
\[\phi = T + a^{q-1}\,T\,\tau + a^{q^2-1}\,\tau^2, \qquad
\psi = T + a^{q-1}\,T^q\,\tau + a^{q^2-1}\,\tau^2, \qquad 
f = 1 + a^{q-1}\,\tau.\]
One readily checks that $f\cdot \phi = \psi\cdot f$. One computes that $\lambda_\nu(\phi) = \deg a$, $\lambda_\nu(\psi)=1+\deg a$, and that $\V^\nu_f(z) = \min\{z, q\, z-\deg a\}$ for all $z\in\R$. We thus have $\V^\nu_f(\lambda_\nu(\phi)) = \deg a$, and so $\lambda_\nu(\psi)>\V^\nu_f(\lambda_\nu(\phi))$.
\end{ex}

\section{Heights of Drinfeld modules}\label{sec:heights}
We now define two notions of heights for Drinfeld modules defined over a finite extension $K$ of $K_0$ and explain their basic properties.

\subsection{Local invariants}\label{sec:local.heights}
Let $\phi$ be a Drinfeld module of rank $r\geqslant1$ defined over $K$. 
For any $a\in A$, we write 
\[ \phi_a = g_{a, 0} + g_{a,1}\,\tau + g_{a,2}\,\tau^2 + \cdots + g_{a,d_a}\,\tau^{d_a}\in K\{\tau\},\]
where $d_a \coloneq r\,\deg a$, and $g_{a,i}\in K$ for all $i$, with $g_{a,0}=a $ and $g_{a, d_a}\neq 0$.

Let $v$ be a place of $K$. 
Equip the completion $F$ of $K$ at $v$ with the (normalized) valuation $\nu = \ord_v$.
We are thus in a position to use the definitions of valuation polygon, maximal break-points, etc. in the previous section: for legibility, we will write $v$ instead of $\nu = \ord_v$ in indices.

For any $a\in A$, one attaches to $\phi_a\in F\{\tau\}$ the following quantity, as in the previous section:  
\begin{equation*} 
\lambda_v(\phi_a) =  -\min\left\{\frac{\ord_v(g_{a,i})-\ord_v(a)}{q^i-1} \,,\ 1\leqslant i\leqslant d_a\right\}\in\Q.
\end{equation*}
It will be convenient to introduce the following variant:
\begin{equation}\label{eq:def.L}
L_v(\phi_a) \coloneq  -\min\left\{\frac{\ord_v(g_{a,i})}{q^i-1} \,,\ 1\leqslant i\leqslant d_a\right\} \in\Q.
\end{equation}
Note that $L_v(\phi_a) = \lambda_v(\phi_a)$ as soon as $\ord_v(a) = 0$.\\

Finally, if $v$ is a \emph{finite} place of $K$, we let 
\begin{equation}\label{eq:def.Lmax}
\L_v(\phi) \coloneq  \max\left\{L_v(\phi_a), \, a\in A\smallsetminus\{0\}\right\}\in\Q.
\end{equation}
This maximum is well-defined.
For any $a\in A$ we have $\phi_{a+1} = 1+ \phi_a$, so that $L_v(\phi_a) = L_v(\phi_{a+1})$: since at least one of $a$, $a+1$ has  valuation $0$ at $v$, we deduce
\begin{align}\label{eq:local.tag.max}
    \max\left\{L_v(\phi_a),\, a\in A\smallsetminus\{0\}\right\} 
    &=\max\left\{L_v(\phi_a),\, a\in A\text{ such that } \ord_v(a)=0\right\} \notag\\
     &=\max\left\{\lambda_v(\phi_a),\, a\in A\text{ such that } \ord_v(a)=0\right\}.
\end{align}

If $\psi = c^{-1}\cdot \phi\cdot c$ for some element $c\in K^\ast$, then one readily checks that
\[\lambda_v(\psi_a) = \lambda_v(\phi_a) - \ord_v(c) \quad \text{ and } \quad
L_v(\psi_a) = L_v(\phi_a) - \ord_v(c),\]
for all $a\in A$. This leads to $\L_v(\psi)=\L_v(\phi)-\ord_v(c)$.
\medskip

We end this subsection by proving a relation between $L_v(\phi_T)$ and $\L_v(\phi)$, which is stated without proof in \cite{Tag93} (see eq. (5.9.2) there).
\begin{lem}\label{lem:comparison.graded.tag}
    Let $\phi$ be a Drinfeld module of rank $r$ over $K$, and let $v$ be a finite place of $K$. 
    Then $L_v(\phi_T)\leqslant \L_v(\phi)$. Moreover,
    assuming  $\phi$ has stable reduction at $v$, we have 
    $\L_v(\phi) = L_v(\phi_T)$. 
\end{lem}
\begin{proof} 
By \eqref{eq:def.Lmax}, it is clear that $L_v(\phi_T)\leqslant \L_v(\phi)$.
We now assume that $\phi$ has stable reduction at $v$, we want to prove that $L_v(\phi_T)\geqslant \L_v(\phi)$, \ie{}, that $L_v(\phi_a)\leqslant L_v(\phi_T)$ for any non zero $a\in A$.
By definition of stable reduction, there exists $c\in K^\ast$ so that that the Drinfeld module $\widetilde{\phi} := c^{-1}\cdot\phi\cdot c$ satisfies the following property: the coefficients $g'_i$ of $\widetilde{\phi}_{\, T} \in K\{\tau\}$  are integral at $v$ ($\ord_v(g'_i)\geqslant 0$) and at least one of them is a unit at $v$ ($\ord_v(g'_i)= 0$ for some $i$) . This property translates as  $L_v(\widetilde{\phi}_{\,T})=0$. 

For any non zero $a\in A$, the coefficients $h_i$ (with $0\leqslant i \leqslant r\,\deg a$) of $\widetilde{\phi}_{\,a}\in K\{\tau\}$ are integral at $v$ (see Lemma 6.1.1 in \cite{Papi}). 
In particular, we have $L_v(\widetilde{\phi}_{\,a})\leqslant 0 = L_v(\widetilde{\phi}_{\,T})$.

Since we have $L_v(\widetilde{\phi}_{\,b})=L_v(\phi_b)-\ord_v(c)$ for any $b\in A$, we have
\[L_v(\phi_a) = L_v(\widetilde{\phi}_{\,a}) + \ord_v(c) \leqslant L_v(\widetilde{\phi}_{\,T}) + \ord_v(c) = L_v(\phi_T),\]
as was to be shown.
\end{proof}

\subsection{Graded height}\label{sec:graded.height}
Let $\phi$ be a Drinfeld $\F_q[T]$-module of rank $r\geqslant1$ defined over $K$. 
Then $\phi$ is characterized by the twisted polynomial
\[ \phi_T = T + g_1\,\tau + g_2\,\tau^2 + \cdots + g_r\,\tau^r\in K\{\tau\},\qquad g_i\in K, \; g_r\neq 0. \]
Let $v$ be a place of $K$. 
We define the \emph{local component at $v$ of the graded height} by 
\begin{equation}\label{eq:def.local.graded.height}
    \hG^v(\phi) \coloneq  L_v(\phi_T) \in\Q.  
\end{equation}
In case $\ord_v(T)=0$, the local term $\hG^v(\phi)$ coincides with $\lambda_v(\phi_T)$. 
For any $c\in K^\ast$, let $\psi = c^{-1}\cdot\phi\cdot c$. A straightforward computation shows that 
\begin{equation}\label{eq:graded.height.isom}
     \hG^v(\psi) = \hG^v(\phi) - \ord_v(c).
\end{equation}

We record the following statement (see \cite[Lemma 6.1.5(1)]{Papi} or \cite[p. 301]{Tag93}):
\begin{prop}\label{prop:graded.height.stable.red} 
Let $\phi$ be as above, and let $v$ be a  place of $K$. 
Then  $\phi$ has stable reduction at $v$ if and only if $\hG^v(\phi)$ is an integer.  
\end{prop}
    \begin{proof} 
    By definition (see section \ref{sec:stablered}), $\phi$ has stable reduction at $v$ if and only if there exists $c\in K^\ast$ such that $L_v(c^{-1}\cdot\phi_T\cdot c) =0$, \ie{}, such that $\ord_v(c)=L_v(\phi_T)$ (see \eqref{eq:graded.height.isom}).
    Since $\ord_v$ is normalized, this happens exactly when $L_v(\phi_T) \in \ord_v(K^\ast) = \Z$.
    \end{proof}

Note that $\hG^v(\phi)= 0$ for all but finitely many places $v\in M_K$. 
We  can then define the \emph{graded height of $\phi$} by 
\begin{equation}\label{eq:def.graded.height}
    \hG(\phi) \coloneq \frac{1}{[K:K_0]}\sum_{v\in M_K} f_v\,\hG^v(\phi).
\end{equation}
By \eqref{eq:graded.height.isom} and the product formula, if $\phi$ and $\psi$ are $K$-isomorphic Drinfeld modules, they have the same graded height. 
We refer to \cite[\S2.2]{BPR21} for more details about the graded height, including a discussion on the link between $\hG(\phi)$ and the $J$-height of $\phi$, \ie{}, the Weil height of the projective point defined by the $j$-invariants of $\phi$.

\subsection{Taguchi height}\label{sec:taguchi.height}
Let $\phi$ be a Drinfeld module of rank $r\geqslant 1$ defined over $K$. We define the Taguchi height of $\phi$ as a sum of local components.

For an infinite place $v\in M_K^\infty$, consider the corresponding field inclusion $\sigma_v:K\hookrightarrow K_v\hookrightarrow \C_\infty$.
Let $\sigma_v\circ\phi : A\to \C_\infty\{\tau\}$ denote the ``base changed'' Drinfeld module. 
By Theorem \ref{theo:corresp.drinfeld.lattices} in section \ref{sec:uniformization}, 
one can unequivocally associate to the Drinfeld module $\sigma_v\circ\phi$ over $\C_\infty$ an $A$-lattice $\Lambda_{\phi,v}\subset \C_\infty$. 
Define the \emph{local component at $v$ of the Taguchi height} by
\begin{equation}\notag
    \hTag^v(\phi) \coloneq  \log_q\left(\covol(\Lambda_{\phi, v})^{-1/r}\right).
\end{equation}
For a finite place $v\in M_K^f$, define the \emph{local component at $v$ of the Taguchi height} to be
\begin{equation}\notag
    \hTag^v(\phi) \coloneq  \left\lceil\L_v(\phi)\right\rceil = \left\lceil \max\left\{L_v(\phi_a), a\in A\smallsetminus\{0\}\right\}\right\rceil \in\Z, 
\end{equation}
where $\lceil\cdot\rceil$ denotes the ceiling function,  
$L_v(\phi_a)$ is given by \eqref{eq:def.L}, and $\L_v(\phi)$ by \eqref{eq:def.Lmax}. 

We then define the \emph{Taguchi height of $\phi$} in the following manner:
\begin{equation}\label{eq:def.taguchi.height}
    \hTag(\phi/K) \coloneq \frac{1}{[K:K_0]}\left(\sum_{v\in M_K^f} f_v\,\hTag^v(\phi)
    +\sum_{v\in M_K^\infty} n_v\,\hTag^v(\phi)
    \right). 
\end{equation}
The sum in \eqref{eq:def.taguchi.height} is actually finite because $\hTag^v(\phi)=0$ for all but a finite number of places $v\in M_K$.
This notion of height was introduced by Taguchi in \cite{Tag91, Tag93}. 
We refer to \cite[Definition~4.1]{Wei17} and \cite[Definition~5.1]{Wei20} 
for more in-depth accounts.

The Taguchi height is sometimes also called the differential height: one can indeed also define $\hTag(\phi)$ as the degree of a certain line bundle of differential forms ``related to'' $\phi$: more precise accounts may be found in \cite[\S5.3]{Tag93} or \cite[D\'efinition 2.8]{DD99}. 
This latter definition draws a clear parallel with the setting of abelian varieties over~$K$, as studied in \cite{GLP26}.

\begin{rem}\label{rem:var.Tag.isom}
    Let $\phi$ be a Drinfeld module as above, and for any $c\in K^\ast$, let $\phi'= c^{-1}\cdot\phi\cdot c$. 
For any place $v$ of $K$, we have
\[\hTag^v(\phi') = \hTag^v(\phi) - \ord_v(c).\]
Indeed, if $v$ is finite, the relation has already been proved in section~\ref{sec:local.heights}. 
If  $v$ is infinite, the lattices $\Lambda_{\phi',v}$, $\Lambda_{\phi,v}$ respectively associated to $\phi'$, $\phi$  satisfy $c\,\Lambda_{\phi',v} = \Lambda_{\phi,v}$, and therefore
$\covol(\Lambda_{\phi,v}) = |c|_v^r\, \covol(\Lambda_{\phi',v})$
by Lemma~\ref{lem:covol.isog}. 

Summing over all places of $K$, we obtain $\hTag(\phi')=\hTag(\phi)$ by the product formula.
This means that the Taguchi height is invariant by $K$-isomorphism.
\end{rem}

\begin{prop}\label{prop:taguchi.height.stable.red} 
Let $\phi$ be as above, and $v$ be a finite  place of $K$. 
If  $\phi$ has stable reduction at $v$, then $\L_v(\phi)$ is an integer, \ie{}, $\hTag^v(\phi) = \L_v(\phi)$.
\end{prop}
    \begin{proof}
        It suffices to combine Proposition \ref{prop:graded.height.stable.red} (which asserts that $L_v(\phi_T)\in\Z$) with Lemma \ref{lem:comparison.graded.tag} (which asserts that $L_v(\phi_T)=\L_v(\phi)$).
    \end{proof}

In what follows, we mainly consider Drinfeld modules $\phi$ with stable reduction everywhere: in that case, by Proposition \ref{prop:taguchi.height.stable.red}, one can rewrite \eqref{eq:def.taguchi.height} as
\begin{equation}
      \hTag(\phi/K) = \frac{1}{[K:K_0]}\left(\sum_{v\in M_K^f} f_v\, \hG^v(\phi)
      + \sum_{v\in M_K^\infty} n_v\,\log_q\left[\covol(\Lambda_{\phi, v})^{-1/r}\right]\right). 
\end{equation}

\subsection{Local comparison of heights}
We gather here a few  comparison results between $\hG^v(\cdot)$ and $\hTag^v(\cdot)$.
\begin{lem}\label{lem:tag.graded.stable}
Let $\phi$ be a Drinfeld module defined over $K$. Let $v$ be a finite place of $K$.
Then $\hG^v(\phi) \leqslant \hTag^v(\phi)$.
Moreover, if $\phi$ has stable reduction at $v$, then $\hG^v(\phi) = \hTag^v(\phi)$.
\end{lem}
    \begin{proof}
    The first assertion follows from the definitions of the local heights and that of $\L_v(\phi)$:
    \[\hG^v(\phi) = L_v(\phi_T) \leqslant \L_v(\phi) \leqslant \lceil \L_v(\phi)\rceil = \hTag^v(\phi).\]
    Assuming $\phi$ has stable reduction at $v$ replaces both inequalities by equalities: see Lemma \ref{lem:comparison.graded.tag} and Proposition \ref{prop:taguchi.height.stable.red}, respectively.
    \end{proof}

\begin{ex} At an infinite place $v$, it appears difficult to find a lower bound on $\hTag^v(\cdot)$ in terms of $\hG^v(\cdot)$.  Let us give an example.

Assume that $K=\F_q(T)$, and that $v=\infty$ is the place at infinity of $K$. We write $K_\infty$ for the completion of $K$ at $v$.
Let $n\geqslant 2$ be an integer. Choose an element $r\in\C_\infty$ such that $|r|_\infty=1$ and such that the vectors $\eta_2:=1$ and $\eta_1:= r\, T^n$ are $K_\infty$-linearly independent.  
Consider the rank $2$ lattice $\Lambda \coloneq A\cdot r\,T^n + A\cdot 1 \subset \C_\infty$  generated by these two vectors. We let $\phi$ denote the Drinfeld module of rank~$2$ defined over $\C_\infty$ associated to $\Lambda$ (see section \ref{sec:uniformization}).  

One checks that the $A$-basis $(\eta_1, \eta_2)$ is a successive minimum basis for $\Lambda$. By definition, we thus have $\hTag^v(\phi) = -\tfrac{1}{2}\,\log_q(|\eta_1|_\infty\,|\eta_2|_\infty)=  -n/2$.

Writing $\phi_T = T + g_1\, \tau + g_2\,\tau^2$ with $g_1, g_2\in\C_\infty$,
Lemma 5.5.7(2) in \cite{Papi} yields
\[ \ord_v(g_2) = q^{n+1}-q^2 -q \qquad \text{and} \qquad \ord_v(g_1) = -q\] 
(keeping in mind that our ordering of successive minimum bases is the inverse of the one in \cite{Papi}).
From there, one easily concludes that $\hG^v(\phi) = q/(q-1)$.    

In this example, $\hG^v(\phi)$ 
is positive and independent of $n$ while, as $n$ grows, 
$\hTag^v(\phi)$ 
tends to~$-\infty$. This obliterates the hope of finding a general lower bound on $\hTag^v(\cdot)$ in terms of $\hG^v(\cdot)$.
\end{ex}

We also give a general version of Lemma 5.1 in \cite{BPR21}: 
\begin{lem}
Let $\phi$ be a Drinfeld module defined over $K$, and let $v$ be a place of $K$.
For any Drinfeld module $\phi'$ defined over $K$ which is $K$-isomorphic to $\phi$, we have
\[ \hG^v(\phi') - \hTag^v(\phi')=\hG^v(\phi) - \hTag^v(\phi).\]
\end{lem}
    \begin{proof} Let $\phi'$ be $K$-isomorphic to $\phi$, we write (as we may) $\phi' = c^{-1}\cdot\phi\cdot c$ for some $c\in K^\ast$. 
    If $v$ is finite, the definitions of the local heights and the discussion in \S\ref{sec:local.heights} yield
    \begin{align*}
        \hG^v(\phi') - \hTag^v(\phi') 
        &= L_v(\phi') - \lceil\L_v(\phi')\rceil 
        = L_v(\phi') -\ord_v(c) -  \lceil\L_v(\phi') - \ord_v(c)\rceil \\
        &= \hG^v(\phi) - \hTag^v(\phi) + \left(\lceil\L_v(\phi)\rceil -\ord_v(c) -  \lceil\L_v(\phi) - \ord_v(c)\rceil\right).
    \end{align*}
    This leads to the desired relation since $\ord_v(c)$ is an integer.
    If $v$ is infinite, it suffices to combine the fact that $\hG^v(\phi') = \hG^v(\phi)-\ord_v(c)$ (see \S\ref{sec:local.heights}) to Remark \ref{rem:var.Tag.isom}. 
    \end{proof}

\subsection{Isogenies and local heights}\label{sec:isog.local.height}
Let $v$ be a finite place of $K$, one may describe the variation of the local Taguchi height at $v$ by isogeny, in terms of the valuation polygon  of the isogeny.
The following is essentially Lemma 2.3 in \cite{Tag93}.

\begin{prop}\label{prop:var.htag.isog} 
Let $v$ be a finite place of $K$.
Let $\phi, \psi$ be two isogenous Drinfeld modules defined over $K$, with stable reduction at $v$.
Let $f:\phi\to\psi$ be a normalized isogeny. Then
\[\hTag^v(\psi) = \V^v_f\big(\hTag^v(\phi)\big).\]
\end{prop}
\begin{proof}
Let $a\in A$ be such that $\ord_v(a)=0$.
Since $f\cdot\phi_a=\psi_a\cdot f$ in $K\{\tau\}$, and since the constant coefficients of $\phi_a, \psi_a$ are units at $v$, Proposition \ref{prop:lambda.poly.isog} yields
$\lambda_v(\psi_a) = \V^v_f(\lambda_v(\phi_a))$. This rewrites as $L_v(\psi_a) = \V^v_f(L_v(\phi_a))$. 
Remembering that $\V^v_f : \R\to\R$ is increasing, we take a maximum over all  $a\in A$  with $\ord_v(a)=0$ (see \eqref{eq:local.tag.max}) to obtain 
$\L_v(\psi) = \V^v_f(\L_v(\phi))$.

The conclusion follows from Proposition \ref{prop:taguchi.height.stable.red}.
\end{proof}

 If one avoids certain places of $K$, one may also describe the variation of the local graded height by isogeny with no stable reduction hypothesis.

\begin{prop}\label{prop:var.hg.isog} 
Let $v$ be a finite place of $K$.
Let $\phi, \psi$ be two isogenous Drinfeld modules defined over $K$.
Let $f:\phi\to\psi$ be a normalized isogeny. Then
\[\hG^v(\psi)  = \V^v_f\big(\hG^v(\phi)\big) .\]
\end{prop}
    \begin{proof}
    Since $v$ is finite, there exists $a\in \{T,T+1\}$ such that $\ord_v(a)=0$.
    For such $a$, the constant coefficient of $\phi_a$ is a unit at $v$, so that $L_v(\phi_a)=\lambda_v(\phi_a)$.
    And since $\phi_a$, $\phi_T$ differ at most by their constant coefficients, we have $L_v(\phi_a) = L_v(\phi_T)$. 
    Hence $\lambda_v(\phi_a) =L_v(\phi_a) = L_v(\phi_T) = \hG^v(\phi)$ and, similarly,
     $\lambda_v(\psi_a) = \hG^v(\psi)$.

    Since $\ord_v(a)=0$, we deduce from Proposition \ref{prop:lambda.poly.isog} that $\lambda_v(\psi_a) = \V^v_f(\lambda_v(\phi_a))$. The previous paragraph allows to rewrite this equality in the desired form.
    \end{proof}

\subsection{Stable reduction and isogenies}
In this subsection, we show that isogenies preserve stable reduction at a place $v\in M_K$. This statement is a general version of \cite[Proposition 2.4]{Tag93}. 
\begin{prop}\label{prop:isog.stable.red}
Let $K$ be a finite extension of $K_0$, and let $\phi, \psi$ be isogenous Drinfeld modules defined over $K$.
For any finite place $v$ of $K$, one has: 

\centerline{$\phi$ has stable reduction at $v$ if and only if $\psi$ has stable reduction at $v$.}
\end{prop}
    \begin{proof}
    Being isogenous is an equivalence relation and is, in particular, symmetric (because of the existence of dual isogenies, \cite[Proposition 3.3.12]{Papi}).
    So it  suffices to prove one of the implications.
    
    Let $f:\phi\to\psi$ be an isogeny defined over $K$ and
    assume that $\phi$ has stable reduction at~$v$. By Proposition~\ref{prop:graded.height.stable.red}, $L_v(\phi_T)=\hG^v(\phi)$ is an integer. 
    We know from the previous Proposition~\ref{prop:var.hg.isog}  that  $\hG^v(\psi) = \V^v_f(\hG^v(\phi))$.
    The definition (see \eqref{eq:def.val.pol}) of $\V^v_f:\R\to\R$, and the fact that $\ord_v$ is normalized by $\ord_v(K^\ast)=\Z$, directly imply that $\V^v_f(\Z)\subseteq\Z$.
    Therefore, $L_v(\psi_T)=\hG^v(\psi)$ is an integer, and 
    Proposition~\ref{prop:graded.height.stable.red} allows to conclude that $\psi$ has stable reduction at $v$.
    \end{proof}

\section{Parallelogram inequalities}
In this section, we gather the results in earlier sections to conclude the proofs of our main theorems.  

\subsection{Parallelogram configurations} \label{sec:parallelogram}
Let $\phi:A\to K\{\tau\}$ be a Drinfeld module of rank $r\geqslant 1$ defined over a finite extension $K$ of $K_0$, and 
let $G, H\in \KERS_K(\phi)$. 
As was already noted in \S\ref{sec:Drinfeld_isogenies}, both $G\cap H$ and $G+H$ are also elements of $\KERS_K(\phi)$.
Proposition \ref{prop:corres.ker.isog} implies the existence and uniqueness of normalized isogenies in $K\{\tau\}$:
\begin{align*}
    \phi &\to\phi/G \ \text{ with kernel } G, &
    \phi &\to\phi/H  \ \text{ with kernel } H, \\
    \phi &\to\phi/(G\cap H)  \ \text{ with kernel } G\cap H, & \text{ and } \quad
    \phi &\to\phi/(G+H)  \ \text{ with kernel } G+H.
\end{align*}
From the inclusions $G\cap H \subset G,  H$  and $G, H \subset G+H$, and from Proposition \ref{prop:facto.isog}, we can write $\phi\to\phi/G$ as the composition $\phi\to\phi/(G\cap H) \to \phi/G$ of two normalized isogenies. 
Similarly, we can write $\phi\to\phi/H$ as the composition $\phi\to\phi/(G\cap H) \to \phi/H$, etc.
The situation is best summarized by the following commutative diagram in which arrows represent normalized isogenies:
\begin{equation} \label{eq:phi_full_parallelogram}
\begin{tikzcd}
&  & \phi \arrow[ddd, bend left =60] \arrow[lldd, bend right=20] \arrow[rrdd, bend left=20] \arrow[d]  &  &   \\
&  & \phi/(G\cap H) \arrow[lld] \arrow[rrd] \arrow[dd]  &  &    \\
\phi/G \arrow[rrd] &  &  &  & \phi/H\,. \arrow[lld] \\
  &  & \phi/(G+H)  &  &                   
\end{tikzcd}
\end{equation}
This leads to a parallelogram of normalized isogenies:
\begin{center}
\begin{tikzcd}
    &  & \phi/(G\cap H) \arrow[lld, "f"'] \arrow[rrd, "g"] \arrow[dd, "\ell"]  &  & \\
    \phi/G \arrow[rrd] &  &  &  & \phi/H \arrow[lld]\,. \\
    &  & \phi/(G+H) &  &                   
\end{tikzcd}
\end{center}
In this parallelogram, note that $\ker(f)\cap\ker(g)$ is trivial: this follows from Lemma \ref{lemm:ker.gcd.lcm} because the isogeny $\phi\to\phi/(G\cap H)$ in diagram \eqref{eq:phi_full_parallelogram} is the GCD in $K\{\tau\}$ of $\phi\to\phi/G$ and $\phi\to\phi/H$.
Thus, the GCD of $f$ and $g$ in $K\{\tau\}$ is $1$, and their LCM in $K\{\tau\}$ is the isogeny $\ell:\phi/(G\cap H) \to \phi/(G+H)$.

If $v$ is a place of $K$ and if $\phi$ has stable reduction at $v$, Proposition \ref{prop:isog.stable.red} ensures that the five Drinfeld modules appearing in diagram \eqref{eq:phi_full_parallelogram}, which are all isogenous to $\phi$, also have stable reduction at $v$.

\subsection{Taguchi height at infinite places}\label{sec:parineq.Tag.infinite}
In this section, we prove:
\begin{theorem}\label{prop:parineq.local.taguchi.infinite}
    Let $\phi$ be a Drinfeld module defined over $K$, and let $G, H\in \KERS_K(\phi)$.
    For any infinite place $v\in M_K^\infty$ of $K$, we have 
    \begin{equation}
        \hTag^v(\phi/(G\cap H)) + \hTag^v(\phi/(G+ H))  = \hTag^v(\phi/G) +\hTag^v(\phi/H).
    \end{equation}
\end{theorem}
    \begin{proof} Let $v$ be an infinite place of $K$ and $\sigma_v:K\into K_v\into\C_\infty$ be the corresponding inclusion of fields. 
    We ``base change'' all the Drinfeld modules (defined over $K$) involved in diagram \eqref{eq:phi_full_parallelogram}  into Drinfeld modules over $\C_\infty$ \textit{via} $\sigma_v$. We still denote them by the same letters.
    By the discussion in~\S\ref{sec:uniformization}, there are $A$-lattices $\Lambda$, $\Lambda_1$, $\Lambda_2$, $\Lambda_\cap$ and $\Lambda_+$ of the same rank $r$ corresponding to each of $\phi$, $\phi/G$, $\phi/H$, $\phi/(G\cap H)$ and $\phi/(G+H)$, respectively.
    The isogenies in diagram \eqref{eq:phi_full_parallelogram} are all normalized. 
    The correspondence between isogenies of Drinfeld modules over $\C_\infty$ and isogenies of $A$-lattices of \S\ref{sec:uniformization} gives rise to a diagram where arrows represent inclusion of lattices:
        
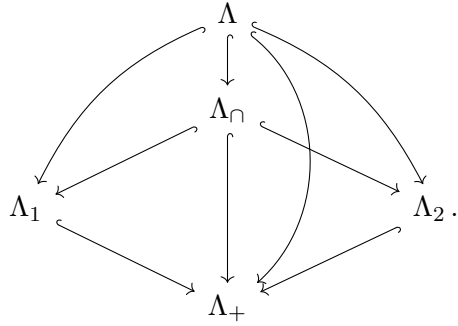
\begin{center}
        \begin{tikzcd}
            &  & \Lambda \arrow[d, hook] \arrow[lldd, hook, bend right=20] \arrow[rrdd, hook, bend left=20] \arrow[ddd, hook, bend left=50] &  & \\
            &  & \Lambda_\cap \arrow[lld, hook] \arrow[rrd, hook] \arrow[dd, hook] &  & \\
            \Lambda_1 \arrow[rrd, hook] &  &  &  & \Lambda_2\,. \arrow[lld, hook] \\
            &  & \Lambda_+ &  &
        \end{tikzcd}
        \end{center}
    Proposition \ref{prop:corresp.isog.lattices} further yields the following four isomorphisms of $\F_q$-vector spaces:
    \[ \Lambda_1/\Lambda \simeq \ker(\phi\to\phi/G) = G,   \quad
    \Lambda_2/\Lambda \simeq \ker(\phi\to\phi/H) = H,  \quad
    \Lambda_\cap/\Lambda \simeq  G\cap H, \quad
    \Lambda_+/\Lambda \simeq  G+H. \]
    From the inclusion $\Lambda_\cap \subseteq\Lambda_1\cap\Lambda_2$ and the isomorphisms $\Lambda_\cap/\Lambda \simeq (\Lambda_1/\Lambda)\cap(\Lambda_2/\Lambda) \simeq (\Lambda_1\cap\Lambda_2)/\Lambda$, 
    we deduce that $\Lambda_\cap = \Lambda_1\cap \Lambda_2$.
    A similar computation yields $\Lambda_+ = \Lambda_1+\Lambda_2$.
    The last displayed diagram simplifies to
    
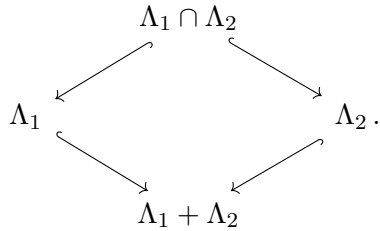
\begin{center}
    \begin{tikzcd}
    & \Lambda_1\cap\Lambda_2 \arrow[ld, hook] \arrow[rd, hook] & \\
    \Lambda_1 \arrow[rd, hook] &    & \Lambda_2 \arrow[ld, hook]\,. \\
    & \Lambda_1+\Lambda_2 &                           
    \end{tikzcd}
    \end{center}
    Applying Theorem \ref{prop:lattice_parallelogram_eq} in this situation yields
    \begin{equation}\label{eq:pre.parineq.Tag.height}
        \log_q \covol(\Lambda_1\cap\Lambda_2)  + \log_q \covol(\Lambda_1 + \Lambda_2) = \log_q \covol(\Lambda_1) + \log_q \covol(\Lambda_2).
    \end{equation}
    By construction of the $A$-lattices here and the definition of the local component at $v$ of the Taguchi height (see \S\ref{sec:taguchi.height}), we have
    \begin{align*} 
    \hTag^v(\phi/G) &= \log_q\left[\covol(\Lambda_1)^{-1/r}\right], \\
    \hTag^v(\phi/H) &= \log_q\left[\covol(\Lambda_2)^{-1/r}\right],\\
    \hTag^v(\phi/(G\cap H)) &= \log_q\left[\covol(\Lambda_\cap)^{-1/r}\right] =\log_q\left[\covol(\Lambda_1\cap \Lambda_2)^{-1/r}\right], \\
    \text{and }\hTag^v(\phi/(G+H)) &= \log_q\left[\covol(\Lambda_+)^{-1/r}\right]=\log_q\left[\covol(\Lambda_1+\Lambda_2)^{-1/r}\right].
    \end{align*}
    Plugging these into \eqref{eq:pre.parineq.Tag.height} concludes the proof.
    \end{proof}

\subsection{Taguchi height at finite places}\label{sec:parineq.local.Tag.finite}
We now turn our attention to finite places of $K$ and prove:
\begin{theorem}\label{prop:parineq.local.taguchi.finite}
    Let $\phi$ be a Drinfeld module defined over $K$, and let $G, H\in \KERS_K(\phi)$.
    For any finite place $v\in M_K$ of $K$ where $\phi$ has stable reduction, we have 
    \begin{equation*}
        \hTag^v(\phi/(G\cap H)) + \hTag^v(\phi/(G+ H))\leqslant\hTag^v(\phi/G) +\hTag^v(\phi/H).
    \end{equation*}
\end{theorem}
\begin{proof}
As discussed in \S\ref{sec:parallelogram}, the data leads to a parallelogram of normalized isogenies
\begin{center}
\begin{tikzcd}
           &  &  \phi/(G\cap H) \arrow[lld, "f"'] \arrow[rrd, "g"] \arrow[dd, "\ell"]&  & \\
 \phi/G \arrow[rrd] &  &                                &  &  \phi/H \arrow[lld] \\
           &  & \phi/(G+H)                             &  &                   
\end{tikzcd}
\end{center}
where $f$ and $g$ have trivial GCD in $K\{\tau\}$, and $\ell \in K\{\tau\}$ is the LCM of $f$ and $g$.

Writing $\mu = \hTag^v(\phi/(G\cap H))$ for simplicity,  Proposition \ref{prop:var.htag.isog} yields
\begin{equation*}
    \hTag^v(\phi/G) = \V^v_f(\mu), \qquad 
    \hTag^v(\phi/H) = \V^v_g(\mu), \quad \text{ and } \quad
    \hTag^v(\phi/(G+H)) = \V^v_\ell(\mu).
\end{equation*}
Applying Theorem \ref{thm:ineq.V.lcm.gcd}, we deduce that $\V^v_\ell(\mu)  \leqslant  \V^v_f(\mu) + \V^v_g(\mu) - \mu$.
In other words, we have
\[ \hTag^v(\phi/(G+H)) \leqslant  \hTag^v(\phi/G)+ \hTag^v(\phi/H)- \hTag^v(\phi/(G\cap H)),\]
which is exactly the desired parallelogram inequality.
\end{proof}

\subsection{Global parallelogram inequality for the Taguchi height}
Summing, with appropriate weights, the parallelogram inequalities for local Taguchi heights obtained in \S\ref{sec:parineq.Tag.infinite} and \S\ref{sec:parineq.local.Tag.finite}, one obtains the parallelogram inequality for the global Taguchi height (Theorem \ref{itheo:main}).
More specifically, Theorem~\ref{itheo:main} is a direct consequence of the definition of the Taguchi height, the equalities in Theorem \ref{prop:parineq.local.taguchi.infinite}, and the inequalities of Corollary \ref{prop:parineq.local.taguchi.finite}.

\subsection{Graded height at a finite place}\label{sec:parineq.local.graded}
 Let $v$ be a finite place of $K$. For any Drinfeld module~$\varphi$ with stable reduction at $v$, our Lemma \ref{lem:tag.graded.stable} ensures that $\hTag^v(\varphi) =\hG^v(\varphi)$. Therefore, one could rewrite Theorem \ref{prop:parineq.local.taguchi.finite} in terms of local graded height. 
However, one can directly prove a parallelogram inequality for the local graded height  (Theorem~\ref{itheo:local.graded}) with no stable reduction assumption, as follows:

\begin{theorem}\label{prop:parineq.local.graded}
Let $\phi$ be a Drinfeld module defined over $K$. Let $v$ be a finite place of $K$. 
For any $G, H\in \KERS_K(\phi)$, we have
    \begin{equation}
        \hG^v(\phi/(G\cap H)) + \hG^v(\phi/(G+ H))\leqslant\hG^v(\phi/G) +\hG^v(\phi/H).
    \end{equation}
\end{theorem}
    \begin{proof}
    The proof is the same as that of Theorem \ref{prop:parineq.local.taguchi.finite}: one simply replaces $\hTag^v(\cdot)$ by $\hG^v(\cdot)$, and one invokes Proposition \ref{prop:var.hg.isog} instead of Proposition \ref{prop:var.htag.isog}.
    \end{proof}

\subsection{Graded height at an infinite place}\label{sec:parineq.graded.infinite}
We  now explain that it is not possible to extend Theorem~\ref{itheo:local.graded} for infinite places by giving an example in which the inequality in Theorem~\ref{itheo:local.graded} fails to hold.
We work at the infinite place $v=\infty$ of $K=K_0=\F_q(t)$ given by $\ord_v(\cdot) = -\deg_T(\cdot)$.
For any Drinfeld module $\phi$ of rank $2$ defined over $K$ by $\phi_T = T +g\,\tau +\Delta \,\tau^2$, we have (see \S\ref{sec:graded.height})
\[\hG^v(\phi) = L_v(\phi_T) = \max\left\{\frac{\deg g}{q-1}, \frac{\deg \Delta}{q^2-1}\right\}.\]

Let $a,b\in\F_q^\ast$, and let $m,n$ be two integers such that $0\leqslant m<n$. 
Consider the two normalized Ore polynomials $f=1-a\,T^m\,\tau$ and $g=1-b\,T^n\,\tau$. 
We write $\delta = b\,T^n-a\,T^m\in K^\ast$.
It is clear that $f,g$ have trivial GCD in $K\{\tau\}$, and one checks that their LCM in $K\{\tau\}$ is $\ell = 1 + x\,\tau + y\,\tau^2$ with
\[ x = \delta^{-q}\,\left(-(b\,T^n)^{q+1}+(a\,T^m)^{q+1}\right), \qquad
y = ab\,\delta^{1-q}\,T^{(m+n)q}.\]
Indeed, we have
$\ell = \left(1 - \delta^{1-q}\,b\,T^{nq}\,\tau\right)\cdot f = \left(1 - \delta^{1-q}\,a\,T^{mq}\,\tau\right)\cdot g$.
We also let
\begin{align*}
 g_0 &= x\,T, &
 g_1 &= -a\,T^{m+q} - b\,\delta^{1-q}\,T^{nq+1},&
 g_2 &= -b\,T^{n+q} - a\,\delta^{1-q}\,T^{mq+1}, &
 g_3 &= x\,T^q, \\
 \Delta_0 &= y\,T, &
 \Delta_1 &= ab\,\delta^{q(1-q)}\,T^{nq^2+m+q},&
 \Delta_2 &= ab\,\delta^{q(1-q)}\,T^{mq^2+n+q},&
 \Delta_3 &= y\,T^q.
\end{align*}
For convenience, we record the degrees of the elements of $K$ introduced so far:

\vspace{-.8em}
\begin{enumerate}[\it(i)]\setlength{\itemsep}{0em}
    \item $\deg x = n$ and $\deg y = n+m\,q$.
    Note that $(q+1)\deg g_0 > \deg\Delta_0$, and $(q+1)\deg g_3 > \deg\Delta_3$.
    \item  $\deg \Delta_2 = m\,q^2+(1+q-q^2)\,n +q $ and $\deg g_2 =  n+q$.
    Note that $(q+1)\deg g_2> \deg \Delta_2$.
    \item $\deg \Delta_1 = m+(n+1)\,q$. 
    If $m\neq n+1-q$, then
    $\displaystyle\deg g_1 = \max\{n+1, m+q\}$, in which case
    we have $(q+1)\deg g_1> \deg\Delta_1$.
    If $m = n+1-q$, the dominant term in $-\delta^q\,g_1 = a\,\delta^q\,T^{n+1} + b\,\delta\,T^{nq+1}$ is $b(a+b)\,T^{nq +n+1}$,
    so that $ \deg g_1 = n+1$ if $a+b\neq 0$, and $ \deg g_1 < n+1$ if $a+b= 0$.
\end{enumerate}

\medskip
With the eight elements $g_i, \Delta_i$ of $K$ introduced above, we define four Drinfeld modules $\phi_0, \phi_1, \phi_2, \phi_3$ of rank $2$ by setting $\phi_{i, T} \coloneq T +g_i \,\tau + \Delta_i \,\tau^2$ for all $i\in\{0, 1, 2, 3\}$.
Tedious but straightforward calculations (the recent preprint \cite{ABCLNP26} provides algorithmic tools which help efficiently carry out these computations)
show that these four Drinfeld modules sit in a parallelogram of isogenies:
\begin{center}
\begin{tikzcd}
           &  &  \phi_0 \arrow[lld, "f"'] \arrow[rrd, "g"] \arrow[dd, "\ell"]&  & \\
 \phi_1 \arrow[rrd] &  &                                &  &  \phi_2 \arrow[lld] \\
           &  & \phi_3                             &  &                   
\end{tikzcd}
\end{center}
All five isogenies are normalized,  $\ker f\cap\ker g$ is trivial, and $\ker \ell = \ker f +\ker g$.
We now compute the local graded heights of these four Drinfeld modules.
By \textit{(i)} above, we have 
\[\hG^v(\phi_0) = \frac{n+1}{q-1} \qquad\text{and}\qquad \hG^v(\phi_3) = \frac{n+q}{q-1} = \hG^v(\phi_0)+1.\]
Similarly, \textit{(ii)} leads to $\displaystyle \hG^v(\phi_2) = \frac{n+q}{q-1} =\hG^v(\phi_3)$.
Finally, we bound $\hG^v(\phi_1)$ using \textit{(iii)}:

\vspace{-.8em}
\begin{itemize}\setlength{\itemsep}{0em}
\item If $m < n+1-q$ (or if $m=n+1-q$ with $a+b\neq 0$), we have
$\hG^v(\phi_1) = \frac{\deg g_1}{q-1} = \hG^v(\phi_0)$.   
\item If $m> n+1-q$, we have
$\hG^v(\phi_1) = \frac{m+q}{q-1} > \frac{n+1}{q-1}  =\hG^v(\phi_0)$.
\item If $m =  n+1-q$ and $a+b = 0$, we have $\deg g_1 < (n+1)$ and $\deg \Delta_1<(q+1)(n+1)$. From which it follows that
$\displaystyle\hG^v(\phi_1) = \max\left\{\frac{\deg g_1}{q-1}, \frac{\deg\Delta_1}{q^2-1}\right\} < \frac{n+1}{q-1} = \hG^v(\phi_0)$.
\end{itemize} 

With these values of $\hG^v(\phi_i)$ at hand, we now observe the following.
If $m \neq n+1-q$ (or if $m=n+1-q$ with $a+b\neq 0$), the parallelogram inequality 
\begin{equation}\label{eq:parineq.graded.example}
    \hG^v(\phi_0) + \hG^v(\phi_3) \leqslant \hG^v(\phi_1) + \hG^v(\phi_2)
\end{equation}
does hold. We remark that \eqref{eq:parineq.graded.example} is actually an \emph{equality} if $m<n+1-q$, and that that the inequality is strict if instead $m>n+1-q$.

If we choose $m=n+1-q$ and $a+b=0$ however, our computations show that
\[\hG^v(\phi_0) + \hG^v(\phi_3) > \hG^v(\phi_1) + \hG^v(\phi_2).\]
In that case, the parallelogram inequality \emph{does not hold}! 
 
\bigskip
 
\centerline{\rule{10cm}{.5pt}}

\vspace{-.8em}
\paragraph{Acknowledgments --} 
The authors are grateful to Francesco Campagna for helpful comments on this paper, and they thank the referee for constructive feedback.
Most of this project was completed during a \textit{GandA Research In Pairs} visit of RG to LB. 
This visit was made possible through funding from IRN GandA (CNRS), which the authors heartily thank.
RG wishes to thank Stellenbosch University for its hospitality during his stay, and CRM-CNRS in Montr\'eal for excellent working conditions during the writing of this paper.
RG is partially supported by ANR-23-CE40-0006 \textit{GAEC}, 
and FP  by ANR-20-CE40-0003 \textit{Jinvariant}.

\centerline{\rule{10cm}{.5pt}}

\small 
\printbibliography{}

@article{BPR21,
 SHORTHAND = {BPR21},
    AUTHOR = {Breuer, Florian and Pazuki, Fabien and Razafinjatovo, Mahefason Heriniaina},
     TITLE = {Heights and isogenies of {D}rinfeld modules},
   JOURNAL = {Acta Arithmetica},
    VOLUME = {197},
      YEAR = {2021},
    NUMBER = {2},
     PAGES = {111--128},
        DOI = {10.4064/aa191029-8-7},
      % ISSN = {0065-1036,1730-6264},
   MRCLASS = {11G09 (11G50 14G17 14G40)},
  MRNUMBER = {4189716},
MRREVIEWER = {Dragos\ Ghioca},
       % URL = {https://doi.org/10.4064/aa191029-8-7},
}

@book {BG06,
SHORTHAND = {BG06},
    AUTHOR = {Bombieri, Enrico and Gubler, Walter},
     TITLE = {Heights in {D}iophantine geometry},
    SERIES = {New Mathematical Monographs},
    VOLUME = {4},
 PUBLISHER = {Cambridge University Press},
      YEAR = {2006},
     %PAGES = {xvi+652},
%      ISBN = {978-0-521-84615-8; 0-521-84615-3},
%   MRCLASS = {11G50 (11-02 11G10 11G30 11J68 14G40)},
%  MRNUMBER = {2216774},
%MRREVIEWER = {Yuri\ Bilu},
       DOI = {10.1017/CBO9780511542879},
       %URL = {https://doi.org/10.1017/CBO9780511542879},
}

@article{DD99,
 SHORTHAND = {DD99},
    AUTHOR = {David, Sinnou and Denis, Laurent},
     TITLE = {Isog\'enie minimale entre modules de {D}rinfel'd},
   JOURNAL = {Mathematische Annalen},
    VOLUME = {315},
      YEAR = {1999},
    NUMBER = {1},
     PAGES = {97--140},
        DOI = {10.1007/s002080050319},
      % ISSN = {0025-5831,1432-1807},
   MRCLASS = {11G09 (11J93)},
  MRNUMBER = {1717545},
MRREVIEWER = {David\ Goss},
       % URL = {https://doi.org/10.1007/s002080050319},
}

@article{Re22,
 SHORTHAND = {R\'em22},
    AUTHOR = {R\'emond, Ga\"el},
     TITLE = {Propri\'et\'es de la hauteur de {F}altings},
   JOURNAL = {Annali della Scuola Normale Superiore di Pisa. Classe di Scienze. Serie V},
    VOLUME = {23},
      YEAR = {2022},
    NUMBER = {3},
     PAGES = {1589--1596},
        DOI = {10.2422/2036-2145.202010_062},
      % ISSN = {0391-173X,2036-2145},
   MRCLASS = {11G10 (11G50)},
  MRNUMBER = {4497754},
MRREVIEWER = {Ulrich\ Derenthal},
       % URL = {https://doi.org/10.2422/2036-2145.202010_062},
}

@article {Tag91,
 SHORTHAND = {Tag91},
    AUTHOR = {Taguchi, Yuichiro},
     TITLE = {Semisimplicity of the {G}alois representations attached to
              {D}rinfel'd modules over fields of ``finite characteristics''},
   JOURNAL = {Duke Math. J.},
  FJOURNAL = {Duke Mathematical Journal},
    VOLUME = {62},
      YEAR = {1991},
    NUMBER = {3},
     PAGES = {593--599},
%      ISSN = {0012-7094,1547-7398},
%   MRCLASS = {11G09 (11R32 11R58)},
%  MRNUMBER = {1104809},
%MRREVIEWER = {Ernst-Ulrich\ Gekeler},
       DOI = {10.1215/S0012-7094-91-06225-3},
%       URL = {https://doi.org/10.1215/S0012-7094-91-06225-3},
}

@article{Tag93,
 SHORTHAND = {Tag93},
    AUTHOR = {Taguchi, Yuichiro},
     TITLE = {Semi-simplicity of the {G}alois representations attached to Drinfel'd modules over fields of ``infinite characteristics''},
   JOURNAL = {Journal of Number Theory},
    VOLUME = {44},
      YEAR = {1993},
    NUMBER = {3},
     PAGES = {292--314},
       DOI = {10.1006/jnth.1993.1055},
      % ISSN = {0022-314X,1096-1658},
   MRCLASS = {11G09 (11R58)},
  MRNUMBER = {1233291},
MRREVIEWER = {David\ Goss},
       % URL = {https://doi.org/10.1006/jnth.1993.1055},
}

@article{Wei20,
 SHORTHAND = {Wei20},
    AUTHOR = {Wei, Fu-Tsun},
     TITLE = {On {K}ronecker terms over global function fields},
   JOURNAL = {Inventiones Mathematicae},
    VOLUME = {220},
      YEAR = {2020},
    NUMBER = {3},
     PAGES = {847--907},
   %MRCLASS = {11G09 (11M36 11R58)},
  %MRNUMBER = {4094971},
%MRREVIEWER = {Gabriel\ D.\ Villa-Salvador},
       DOI = {10.1007/s00222-019-00944-8},
       % URL = {https://doi.org/10.1007/s00222-019-00944-8},
}

@article{GLP26,
  shorthand     = {GLP26},
  title         = {Variation of height in an isogeny class over a function field},
  author        = {Griffon, Richard and Le Fourn, Samuel and Pazuki, Fabien},
  journal       = {Documenta Math, to appear},
  year          = {2026},
  eprint        = {2503.14318},
  archivePrefix = {arXiv},
  primaryClass  = {math.NT},
  DOI = {10.4171/dm/1094}
  % url       = {https://arxiv.org/abs/2503.14318}, 
}

@book{Papi,
 SHORTHAND = {Pap23},
    AUTHOR = {Papikian, Mihran},
     TITLE = {Drinfeld modules},
    SERIES = {Graduate Texts in Mathematics},
    VOLUME = {296},
 PUBLISHER = {Springer},
      YEAR = {2023},
  %    ISBN = {978-3-031-19706-2},
 %  MRCLASS = {11G09 (11-02 11R58)},
 % MRNUMBER = {4592575},
%MRREVIEWER = {Ernst-Ulrich\ Gekeler},
 DOI = {10.1007/978-3-031-19707-9},
       % URL = {https://doi.org/10.1007/978-3-031-19707-9},
}

@book{Gos12,
  shorthand = {Gos12},
  title     = {Basic {S}tructures of {F}unction {F}ield {A}rithmetic},
  author    = {Goss, David},
  year      = {2012},
  publisher = {Springer},
 doi       = {10.1007/978-3-642-61480-4},
  %isbn      = {978-3-540-63541-3},
  %pagetotal = {424}
}

@book{Lang,
  shorthand = {Lang},
  title     = {Algebra},
  author    = {Lang, Serge},
  series    = {Graduate Texts in Mathematics},
  volume    = {211},
%  pagetotal = {914},
  year      = {2002},
  publisher = {Springer},
%  %isbn      = {978-0-387-95385-4},
 doi       = {10.1007/978-1-4613-0041-0}
}

@article{Pazuki,
  shorthand = {Paz19},
    AUTHOR = {Pazuki, Fabien},
     TITLE = {Modular invariants and isogenies},
   JOURNAL = {International Journal of Number Theory},
    VOLUME = {15},
      YEAR = {2019},
    NUMBER = {3},
     PAGES = {569--584},
  % MRCLASS = {11G50 (11G05 14G40 14J15 14K02)},
 % MRNUMBER = {3925753},
%MRREVIEWER = {Yuri\ Bilu},
   DOI = {10.1142/S1793042119500295},
%        URL = {https://doi-org.ezproxy.uca.fr/10.1142/S1793042119500295},
}

@article{Wei17,
  shorthand = {Wei17},
%  ISSN = {00029327, 10806377},
  author = {Wei, Fu-Tsun},
  journal = {American Journal of Mathematics},
  number = {4},
  pages = {1047--1084},
  publisher = {The Johns Hopkins University Press},
  title = {Kronecker limit formula over global function fields},
  volume = {139},
  doi = {10.1353/ajm.2017.0027},
  year = {2017}
}

@incollection{Robba,
    shorthand = {Rob85},
    AUTHOR = {Robba, Philippe},
     TITLE = {Lemmes de {H}ensel pour les op\'erateurs diff\'erentiels et
              les op\'erateurs aux diff\'erences (d'apr\`es {A}. {D}uval)},
 BOOKTITLE = {Algebra colloquium (Publications de l'Institut de recherche math{\'e}matiques de Rennes)},
     PAGES = {10--24},
 PUBLISHER = {Univ. Rennes I, Rennes},
      YEAR = {1985},
       %URL = {https://www.numdam.org/item/PSMIR_1985___4_10_0.pdf}
}

@book {Rosen,
 shorthand = {Ros02},
    AUTHOR = {Rosen, Michael},
     TITLE = {Number theory in function fields},
    SERIES = {Graduate Texts in Mathematics},
    VOLUME = {210},
 PUBLISHER = {Springer-Verlag, New York},
      YEAR = {2002},
%     PAGES = {xii+358},
      % ISBN = {0-387-95335-3},
   % MRCLASS = {11R58 (11R60 11T55)},
  % MRNUMBER = {1876657},
% MRREVIEWER = {Ernst-Ulrich\ Gekeler},
       DOI = {10.1007/978-1-4757-6046-0},
%       URL = {https://doi-org.ezproxy.uca.fr/10.1007/978-1-4757-6046-0},
}

@incollection {Hayes,
shorthand = {Hay79},
    AUTHOR = {Hayes, David R.},
     TITLE = {Explicit class field theory in global function fields},
 BOOKTITLE = {Studies in algebra and number theory},
    SERIES = {Adv. Math. Suppl. Stud.},
    VOLUME = {6},
     PAGES = {173--217},
 PUBLISHER = {Academic Press, New York-London},
      YEAR = {1979},
      %ISBN = {0-12-599153-3},
   %MRCLASS = {12A90 (10D30 12A65)},
  %MRNUMBER = {535766},
%MRREVIEWER = {W.-D.\ Geyer},
}

@article{ABCLNP26,
shorthand = {ABC+26},
  title={A computational approach to Drinfeld modules},
  author={Armana, C{\'e}cile and Berardini, Elena and Caruso, Xavier and Leudi{\`e}re, Antoine and Nardi, Jade and Pazuki, Fabien},
  journal={ArXiv preprint ArXiv:2601.02162},
  year={2026},
  url       = {https://arxiv.org/pdf/2601.02162}, 
}

	\normalsize\vfill
	\noindent\rule{7cm}{0.5pt}
 
    \smallskip
	\noindent
	{Liam {\sc Baker}} {(\it \href{liambaker@sun.ac.za}{liambaker@sun.ac.za})}  --
    {\sc Department of Mathematical Sciences, Stellenbosch University}, Stellenbosch, South Africa 7599.
	
    \medskip
	\noindent
	{Richard {\sc Griffon}} {(\it \href{richard.griffon@uca.fr}{richard.griffon@uca.fr})}  --
	{\sc Laboratoire de Math\'ematiques Blaise Pascal, Universit\'e Clermont Auvergne}. 
	3 Place Vasarely, 63710 Aubi\`ere (France).
	
	\medskip
	\noindent
	{Fabien {\sc Pazuki}} {(\it \href{fpazuki@math.ku.dk}{fpazuki@math.ku.dk})} --
	{\sc Department of Mathematical Sciences, University of Copenhagen}.
	Universitetsparken 5, 2100 Copenhagen \o{} (Denmark).

\end{document}